\documentclass{amsart} 
\newcommand{\f}{\frac}

\newcommand{\del}{\partial}

\newcommand{\R}{\mathbb R}

\newcommand{\eps}{\varepsilon}
\renewcommand{\epsilon}{\varepsilon}

\newcommand{\Om}{\Omega}

\newcommand{\Jeps}{U_{\eps}}

\newcommand{\dist}{\operatorname{dist}}
\newcommand{\Oeps}{\Omega_\eps^{\mathrm{p}}}
\newcommand{\Heps}{\mathcal{H}_{\eps}}
\newcommand{\U}{\mathcal{U}_{\eps}}
\newcommand{\V}{\mathcal{V}_{\eps}}
\newcommand{\Vt}{\tilde{\mathcal{U}}_{\eps}}
\newcommand{\W}{\mathcal{W}_{\eps}}
\newcommand{\diam}{\operatorname{diam}}

\newcommand{\xmax}{9}
\newcommand{\boxy}{-1.8}

\usepackage[english]{babel}
\usepackage{tikz}
\usepackage{pgfplots}
\usetikzlibrary{arrows,decorations.pathmorphing,backgrounds,positioning,fit,petri,patterns}
\usepackage[utf8]{inputenc}
\usepackage{graphicx}
\usepackage[colorlinks=true,linkcolor=black,citecolor=black]{hyperref}
\usepackage{amsfonts}
\usepackage{enumerate}
\usepackage{booktabs}
\usepackage{array} 
\usepackage{paralist} 
\usepackage{subfig} 
\usepackage{mathtools}
\usepackage{tabu}
\usepackage{amsthm}
\usepackage{empheq}
\usepackage{amsopn}
\usepackage{dsfont}
\usepackage{wrapfig}
\usepackage[font=small]{caption}
\usepackage{units}
\usepackage[symbol]{footmisc}
\usepackage{todonotes}
\usepackage{amssymb}
\pgfplotsset{compat=1.10}
\usepgfplotslibrary{fillbetween}

\allowdisplaybreaks

\theoremstyle{definition}
\newtheorem{de}{Definition}[section]

\theoremstyle{plain}
\newtheorem{prop}[de]{Proposition}
\newtheorem{lemma}[de]{Lemma}
\newtheorem{theorem}[de]{Theorem}
\newtheorem{corollary}[de]{Corollary}

\numberwithin{equation}{section}
\theoremstyle{remark}
\newtheorem{remark}[de]{Remark}

\title{\sc A Strange Vertex Condition Coming from Nowhere}
\author{F. R\"osler}
\email{RoslerF@cardiff.ac.uk}
\address{School of Mathematics, Cardiff University, Senghennydd Road, Cardiff CF24 4AG, Wales, UK}
\thanks{
	The author thanks P. Dondl and K. Cherednichenko for discussions and inspiration.
	Furthermore, the author would like to thank the anonymous reviewers and the editors for helpful suggestions. Finally the author acknowledges support from the European Union's Horizon 2020 Research and Innovation Programme under the Marie Sk{\l}odowska-Curie grant agreement No. 885904.
}

\begin{document}

\maketitle 
\begin{abstract}
	We prove norm-resolvent and spectral convergence in $L^2$ of solutions to the Neumann Poisson problem $-\Delta u_\eps = f$ on a domain $\Om_\eps$ perforated by Dirichlet-holes and shrinking to a 1-dimensional interval. The limit $u$ satisfies an equation of the type $-u''+\mu u = f$ on the interval $(0,1)$, where $\mu$ is a positive constant.
	
	As an application we study the convergence of solutions in perforated graph-like domains. We show that if the scaling between the edge neighbourhood and the vertex neighbourhood is chosen correctly, the constant $\mu$ will appear in the vertex condition of the limit problem. In particular, this implies that the spectrum of the resulting quantum graph is altered in a controlled way by the perforation.
\end{abstract}
\section{Introduction} 
Let $N\geq 3$ and consider an open subset $\Om_\eps$ of $\R^N$ of the form $\Om_\eps = \eps\Om_0\times(0,1)$ {(see Section \ref{sec:geometric_setting} for precise definitions)}. Let us introduce a perforation of this domain by removing periodically distributed spherical holes of distance {$\delta_\eps\in(0,\eps)$} (cf. Figure \ref{fig:thin_perforated_domain}). On this domain we consider the Poisson equation with Dirichlet boundary conditions on the holes of radius {$r_\eps\ll \delta_\eps$}. We ask the question whether the solutions $u_\eps$ to this equation converge in a meaningful sense to a function $u$ on the interval $(0,1)$ and whether $u$ is the solution of a reasonable ``limit'' differential equation. 

Homogenisation problems of a similar type have been studied extensively for a long time \cite{CM,RT,MK} and recently gained more attention 
	(cf. \cite{Jikov2000,Pas06} for perforated domains of fixed size with Neumann boundary conditions,  \cite{MS10} for perforated domains with periodic boundary conditions, \cite{BCD} for domains perforated along a curve. Advances towards operator norm  and spectral convergence in perforated domains have been made in \cite{Pas06,BCD,CDR,KP}).
	 A result by Cioranescu \& Murat gives a positive answer to the question of convergence of solutions in the case where the size of $\Om_\eps$ remains constant, but the holes shrink and concentrate. In fact, they showed that the solutions of $-\Delta u_\eps = f$ converge strongly in $L^2(\Om)$ to the solution $u\in H^1_0(\Om)$ of $(-\Delta+\bar\mu)u=f$, where $\bar\mu>0$ is a constant related to the harmonic capacity of the unit ball. The constant $\mu$ (which was dubbed a ``strange term coming from nowhere'' in \cite{CM}) will appear frequently in later sections of this article and we will henceforth refer to $\mu$ as the \emph{strange term}.

	The general idea of coupling thin geometry with a highly oscillating boundary of the domain has also gained interest during the last decade. Indeed, elliptic problems on a thin domain whose boundary is given as the graph of a rapidly oscillating function $G_\eps$ have been studied in \cite{AP10,AV14,AV16}.
	The more specific situation of a perforated thin domain was the object of study in \cite{MP10,MP12} (see also the references therein). The effects of perforations in thin domains on spectral gaps have been studied in \cite{N10}. 
	
	The present article differs from these works in several ways. First, the geometric situation is different in the sense that the radius of the holes does not have the same scaling as the distance between the holes or the thickness of the domain. Second, the boundary conditions we consider on the surface of the holes are Dirichlet (rather than Neumann), which changes the analysis of the problem completely and ultimately leads to the appearance of the strange term $\mu$ in the limiting equation. Moreover, the emphasis of the present work differs from those mentioned in the last two paragraphs. We take an operator theoretic point of view and prove that the operators involved converge in \emph{norm-resolvent sense}, i.e. the resolvents of the operator family indexed by $\eps$ converge in the uniform operator topology. This notion of convergence is stronger than that of \emph{strong convergence}, which is more commonly studied in classical homogenisation theory. In particular, norm-resolvent convergence implies a number of physically interesting consequences like local convergence of spectra (cf. Section \ref{SpSec}) or convergence of the associated semigroups. Finally, our results are applied to so-called \emph{graph-like domains} in Section \ref{sec:graph-like_Domains}, where the additional challenge of determining vertex conditions for the limiting equation is present. This situation is similar to that in \cite{P06}, however, there the author did not consider the effect of perforations.
	
	This article is organised as follows. In Section \ref{sec:geometric_setting}, we give a precise description of the geometric situation at hand and the resulting boundary value problem in the perforated thin domain. Section \ref{sec:main_results} contains the statements of our main theorems and relevant corollaries. Sections \ref{sec:general_convergence_results}, \ref{sec:proof_of_theorem} and \ref{NormandSpec} are devoted to the proof of our main theorem. In Section \ref{SpSec} we prove local convergence of spectra as a corollary of norm-resolvent convergence. Finally, in Section \ref{sec:graph-like_Domains} we apply our results to perforated graph-like domains and obtain vertex conditions for the limiting problem on the underlying metric graph.

\section{Geometric setting}\label{sec:geometric_setting}

In this article we consider the following homogenisation problem. Let $N\geq 3$ and $\Om_0\subset\R^{N-1}$ be a bounded open set with $\del\Om_0$ of class $C^2$ and let $\Om:=\Om_0\times(0,1)$. For $\eps>0$, let $\delta_\eps < \eps$ and define the set $\tilde T_\eps:=\bigcup_{i\in 2\delta_\eps\mathbb Z^N}B_{r_\eps}(i)$, where $r_\eps = \delta_\eps^{\nicefrac{N}{(N-2)}}$. We consider the domain $\Om_\eps:=\eps\Om_0\times(0,1)$, perforated by the  $B_{r_\eps}(i)$ and shrinking towards a thin rod {as $\eps\to 0$}. 
\begin{figure}[htbp]
	\centering
%

\begin{tikzpicture}[>=stealth]
	\clip (-2,-3) rectangle (10,3);
	
	\shadedraw[top color=gray!20, bottom color=gray!70, xshift=-200] plot [smooth cycle, tension=1] coordinates {(7.4,-0.5)(7.2,0.5) (7,2) (6.7,0) (7,-2)};

	\shade[top color=gray!20, bottom color=gray!70] (0,-2.01) rectangle (8,2.01);
	\draw (0,2.01) -- (8,2.01);
	\draw (0,-2.01) -- (8,-2.01);

	\filldraw[fill = gray!60, xshift=28] plot [smooth cycle, tension=1] coordinates {(7.4,-0.5)(7.2,0.5) (7,2) (6.7,0) (7,-2)};
		
	\draw (0,2.5)node[]{${\Omega_\varepsilon^{\mathrm p}}$} ;
	
	\foreach \y in {-1,1,0} {
	\foreach \x in {1,2,3,4,5,6,7} {
		\draw[xshift=-10, yshift=-7, opacity=0.7] (\x,\y) circle (1mm);
		\shade[ball color=white, xshift=-10, yshift=-7, opacity=0.7] (\x,\y) circle (1mm);
	}}
	
	\foreach \x in {(1,0),(2,0),(3,0),(4,0),(5,0),(6,0),(1,1),(2,1),(3,1),(4,1),(5,1),(6,1),(1,-1),(2,-1),(3,-1),(4,-1),(5,-1),(6,-1),(7,0),(7,1),(7,-1)} {
		\draw[xshift=-3, yshift=0, opacity=0.5] \x circle (1mm);
		\shade[ball color=white, xshift=-3, yshift=0, opacity=0.5] \x circle (1mm);
		}
	
	\foreach \x in {(1,0),(2,0),(3,0),(4,0),(5,0),(6,0),(1,1),(2,1),(3,1),(4,1),(5,1),(6,1),(1,-1),(2,-1),(3,-1),(4,-1),(5,-1),(6,-1),(7,0),(7,1),(7,-1)} {
		\draw[xshift=4, yshift=7, opacity=0.1] \x circle (1mm);
		\shade[ball color=white, xshift=4, yshift=7, opacity=0.4] \x circle (1mm);
		}

	\draw (8,0) -- (9,1);
	\draw (9.4,1.2) node[]{$\varepsilon\Omega_0$};
	
	\draw (-1,0) node[]{$\sim\varepsilon$};
	\draw[->] (-1,0.2) -- (-1,2);
	\draw[->] (-1,-0.2) -- (-1,-2);
	
	\draw[dotted,xshift=-4] (3.04,-1) -- (3.04,-2.6);
	\draw[dotted,xshift=-4] (4.04,-1) -- (4.04,-2.6);
	\draw[xshift=-4] (3.54,-2.4) node[]{$2\delta_\varepsilon$};
	\draw[->,xshift=-4] (2.54,-2.4) -- (3.04,-2.4);
	\draw[->,xshift=-4] (4.54,-2.4) -- (4.04,-2.4);

	\draw[xshift=-4] (6.05,-0.95) -- (6.5,-2.3);
	\draw[xshift=-4] (6.6,-2.6) node[]{$B_{r_\varepsilon}(i)$};
	
\end{tikzpicture}
	\caption{A sketch of the thin perforated domain in 3d.}\label{fig:thin_perforated_domain}
\end{figure}
To this end, define the subset of lattice points which are sufficiently far from the boundary $L_\eps:=\{i\in 2\delta_\eps\mathbb Z^N : \dist(i,\del(\Om_\eps))>\delta_\eps\}$ and the corresponding ``holes'' $T_\eps:=\bigcup_{i\in L_\eps}B_{r_\eps}(i)$. Finally, define the perforated domain
\begin{align*}
	\Om_\eps^\text{p} := \Om_\eps\setminus T_\eps.
\end{align*}
In order to compare functions defined on different domains $\Om_\eps$ and $(0,1)$ we define the operator family
\begin{align*}
	U_\eps:L^1((0,1))&\to L^1(\Om_\eps)\\
	U_\eps \phi &= |\eps\Om_0|^{-\f12}\phi^*,
\end{align*}
where $\phi^*$ denotes the extension of $\phi$ to a constant on every slice $\{t\}\times\eps \Om_0$. Restrictions of $U_\eps$ to subspaces of $L^1(\Om_\eps)$ will also be denoted $U_\eps$. Note that the scaling {$|\eps\Om_0|^{-\f12}$} in the definition of $U_\eps$ was chosen such that for $\phi\in L^2((0,1))$ the norm
	$\|U_\eps\phi\|_{L^2(\Om_\eps)}$ is of order 1 as $\eps\to 0$.
On the domain $\Om_\eps^\mathrm{p}$ we consider the following problem
\begin{align}\label{pde}
\begin{cases}
	\hfill(-\Delta+z) u_\eps = f_\eps,\phantom{0} &\text{in }\Oeps \\
	\hfill u_\eps = 0,\phantom{f_\eps} &\text{on }\del T_\eps \\
	\hfill\del_\nu u_\eps = 0,\phantom{f_\eps} &\text{on }\del\Om_\eps,
\end{cases}
\end{align}
where $z>0$ and $f_\eps\in L^2(\Om_\eps)$ is a family such that $\|f_\eps-\Jeps f\|_{L^2(\Om_\eps)}\to 0$ for some $f\in L^2((0,1))$. This problem can easily be seen to possess a unique solution for each fixed $\eps>0$ by virtue of the Lax-Milgram theorem.

Moreover, let $\Heps:=H^1(\Om_\eps)$ and 
$$
	\mathcal H_\eps^0:=\overline{\bigl\{\phi|_{\Om_\eps} : \phi\in C^\infty_0\big(\R^N\setminus T_\eps\big)  \bigr\}},
$$
where the closure is taken in the $H^1(\Om_\eps)$-norm (this is the space of functions vanishing on the holes). 
For a function $u\in\mathcal H_\eps^0$ we will not distinguish in notation between $u$ and its extension by zero to $\Om_\eps$ (which {belongs to} $\Heps$).

Finally, the following notation will be used frequently. For $x\in\Om_\eps$ we write $x=(\bar x,x_N)$, where $\bar x\in \eps\Om_0$ and $x_N\in(0,1)$. Accordingly, we denote by $\bar\nabla$ the gradient {with respect to} $\bar x$ and by $\del_N$ the partial derivative {with respect to} $x_N$. The {transversally} constant extension of a function $\phi$ from $(0,1)$ to $\Om_\eps$ will be denoted $\phi^*(\bar x,x_N):=\phi(x_N)$. A variable in $(0,1)$ will often be denoted by $t$.

\section{Main results}\label{sec:main_results}

In the above setting, we are going to prove the following results

\begin{theorem}\label{mainth}
	The solutions $u_\eps$ of \eqref{pde} converge to a function $u\in H^1((0,1))$ in the sense that
	$$
		\left\|u_\eps-\Jeps u\right\|_{L^2(\Om_\eps)} \to 0,
	$$
	as $\eps\to 0$ and $u$ solves the ordinary differential equation
	\begin{align}\label{limeq}
	\begin{cases}
		\left(-\f{d^2}{dt^2} + z + \mu\right) u = f, & \text{ in }(0,1)\\
		\hfill u'=0, & \text{ on }\del(0,1),
	\end{cases}
	\end{align}
	where $\mu=2^{-N}S_N(N-2)$, $S_N$ being the surface area of the unit sphere in $\R^N$. 
\end{theorem}
The above theorem can be understood as strong operator convergence $-\Delta_{\Oeps}\xrightarrow{s}-\f{d^2}{dt^2}+\mu$. The next result shows that even a stronger type of convergence holds.
\begin{theorem}
	The above convergence even holds in the norm-resolvent sense.
\end{theorem}
The meaning of ``convergence in the norm-resolvent sense'' will be made precise in Section \ref{NormandSpec} (see Theorem \ref{normconv}). An important corollary of norm-resolvent convergence is convergence of spectra.
\begin{corollary}[Spectral Convergence]\label{SpCon}
	Choose $z=1$ and let $\lambda_k^\eps$ and $\lambda_k$ denote the $k$-th eigenvalues of problem \eqref{pde} and \eqref{limeq}, respectively. There exist a constant $C>0$ and a function $a(\eps)$ with $a(\eps)\to 0$ as $\eps\to 0$ such that 
	\begin{align*}
		|(\lambda_k^\eps)^{-1} - \lambda_k^{-1}| \leq C a(\eps)\qquad\text{ for all }k\in\mathbb N,
	\end{align*}
	where $C$ is independent of $\eps$ and $k$.
\end{corollary}
This corollary will be proved in Section \ref{SpSec}.
The appearance of the additive term $\mu u$ in \eqref{limeq} has been first observed in the classical situation of a perforated domain $\Om$ of fixed size by \cite{MK,CM} and has been dubbed a ``strange term coming from nowhere''. We will in the following refer to $\mu$ as the \emph{strange term}.
\paragraph{Graph-like Domains.}
The above results will be applied to graph-like domains in Section \ref{sec:graph-like_Domains}. In particular, we will show that for a graph-like domain in which the volumes of the fattened edges and the fattened vertices have the same scaling as $\eps\to 0$, the limit will be a quantum graph with vertex conditions of Robin type with parameter $\mu$. For details, see Section \ref{sec:borderline}.

\section{General convergence results on $\Om_\eps$}\label{sec:general_convergence_results}
In the following sections we will prove Theorem \ref{mainth}. We start with some general lemmas about convergence in shrinking domains.
\begin{de}\label{def:strong_convergence}
	A sequence $\phi_\eps\in\Heps$ is said to \emph{strongly converge} to $\phi\in H^1((0,1))$ (we write $\phi_\eps\xrightarrow{H^1}\phi$), if
	\begin{align*}
		 \|\phi_\eps-U_\eps\phi\|_{L^2(\Om_\eps)}^2 + \eps^2 \|\bar \nabla\phi_\eps-\bar\nabla U_\eps\phi\|_{L^2(\Om_\eps)}^2 + \|\del_N\phi_\eps-\del_N U_\eps\phi\|_{L^2(\Om_\eps)}^2 \to 0
	\end{align*}
	as $\eps\to 0$. Strong convergence in $L^2$ is defined analogously, {for which we will write $\phi_\eps\xrightarrow{L^2}\phi$}.
\end{de}
\begin{de}\label{weakdef}
	A sequence $u_\eps\in \Heps$ is said to be \emph{weakly convergent} in $H^1$ to $u\in H^1((0,1))$ (we write $u_\eps\xrightharpoonup{H^1} u$), if for all $\phi_\eps\in\Heps$ with $\phi_\eps\xrightarrow{H^1}\phi$ one has
	\begin{align*}
		\langle u_\eps , \phi_\eps \rangle_{L^2(\Om_\eps)}  + \eps^2\langle \bar\nabla u_\eps , \bar\nabla\phi_\eps \rangle_{L^2(\Om_\eps)} + \langle \del_N u_\eps , \del_N \phi_\eps \rangle_{L^2(\Om_\eps)} \to \langle u , \phi \rangle_{H^1((0,1))}.
	\end{align*}
	Weak convergence in $L^2$ is defined analogously, {for which we will write $\phi_\eps\xrightharpoonup{L^2}\phi$}.
\end{de}
It can easily be seen that in the above sense strong convergence implies weak convergence. 
\begin{remark}
\begin{enumerate}[(i)]
	\item We remark that the concepts of convergence introduced in Definitions \ref{def:strong_convergence} and \ref{weakdef} are not new. Indeed, convergence of sequences in varying Banach spaces has been studied for several decades and Definitions \ref{def:strong_convergence} and \ref{weakdef} are special cases of what is known as \emph{discrete convegrence} (cf. \cite{Stummel1}). Properties of discretely converging sequences of vectors have been studied in the classical works \cite{Stummel1,Stummel2,V81}. In fact, Proposition \ref{weakconv} (i) below is a consequence of \cite[Prop. 1.5]{V81}. We nevertheless chose to include these definitions and proofs in our article in order to keep the presentation as clear and {self-contained} as possible.
	\item The convergence of \emph{operators} defined on varying spaces has also been studied in \cite{Stummel1,Stummel2,V81} to a certain extent. Classical results include various conditions for the strong discrete convergence of bounded operators (and strengthened versions thereof). Let us stress again that in our situation we are dealing with \emph{unbounded} operators for which we are studying the stronger notion of \emph{operator norm} convergence.
	For more recent results on the convergence (especially spectral convergence) of unbounded operators on varying Hilbert spaces, the interested reader may consult \cite{P06,MNP} and \cite{B17,B18}.
\end{enumerate} 
\end{remark}
The next proposition shows that compact embeddings also generalise to shrinking domains.
\begin{prop}\label{weakconv}
	Let $u_\eps\in\Heps$ be a sequence and let there exist a $C>0$ such that 
	\begin{align}\label{H1bdd}
		 \|u_\eps\|_{L^2(\Om_\eps)}^2 + \eps^2 \|\bar \nabla u_\eps\|_{L^2(\Om_\eps)}^2  + \|\del_N u_\eps\|_{L^2(\Om_\eps)}^2 \leq C.
	\end{align}
	for all $\eps>0$. Then 
	\begin{enumerate}[(i)]
	\item there exists a subsequence (still denoted by $u_\eps$) such that $u_\eps \xrightharpoonup{H^1} u$ for some $u\in H^1((0,1))$;
	\item if in addition $\eps^{2}\|\bar \nabla u_\eps\|_{L^2(\Om_\eps)}^2 \to 0$, then one has $\left\|u_\eps-U_\eps u\right\|_{L^2(\Om_\eps)}\to 0$.
	\end{enumerate}
\end{prop}
\begin{proof}
	We use scaling in order to keep the domain fixed. Let $\tilde u_\eps: \Om\to\R,\;  \tilde u_\eps(x):=u_\eps(\eps\bar x,x_N)$. By the usual dilation formula and chain rule we find
	\begin{align*}
		\|u_\eps\|_{L^2(\Om_\eps)}^2 &= \eps^{N-1}\|\tilde u_\eps\|_{L^2(\Om)}^2\\
		\|\del_N u_\eps\|_{L^2(\Om_\eps)}^2 &= \eps^{N-1}\|\del_N \tilde u_\eps\|_{L^2(\Om)}^2\\
		\|\bar \nabla u_\eps\|_{L^2(\Om_\eps)}^2 &= \eps^{N-3}\|\bar\nabla\tilde u_\eps\|_{L^2(\Om)}^2.
	\end{align*}
	Our assumption \eqref{H1bdd} immediately yields $\eps^{N-1}\|\tilde u_\eps\|_{H^1(\Om)}^2\leq C$. Thus, there exists a subsequence $\eps^{\f{N-1}{2}}\tilde u_\eps\rightharpoonup \tilde u$ in $H^1(\Om)$ (in the usual sense). 
	
	Now let $\phi_\eps\in\Heps$ with $\phi_\eps\xrightarrow{H^1}\phi\in H^1((0,1))$. By scaling arguments similar to the above, one immediately obtains that denoting $\tilde\phi_\eps(x):=\phi_\eps(\eps\bar x, x_N)$ and $\phi^*(x):=\phi(x_N)$ one has
	\begin{align*}
		\eps^{\f{N-1}{2}}\tilde\phi_\eps \to \phi^*\quad\text{ strongly in }H^1(\Om).
	\end{align*}
	Consequently,
	\begin{align*}
		\eps^{N-1}\langle \tilde u_\eps , \tilde\phi_\eps\rangle_{H^1(\Om)} \to \langle\tilde u,\phi^*\rangle_{H^1(\Om)}.
	\end{align*}
	Undoing the scaling this can be written as
	\begin{align}\label{welldef}
		\langle u_\eps,\phi_\eps\rangle_{\!L^2(\Om_\eps)} + \eps^{2}\langle \bar\nabla u_\eps,\bar\nabla\phi_\eps\rangle_{\!L^2(\Om_\eps)} 
		 + \langle \del_N u_\eps,\del_N\phi_\eps\rangle_{\!L^2(\Om_\eps)} 
		\;&\to\; 
		\langle \tilde u,\phi^*\rangle_{\!H^1(\Om)}\\
		&= \left\langle \int_\Om \tilde u(\overline{x},\cdot)\,d\overline{x}\, ,\, \phi \right\rangle_{H^1((0,1))},
	\end{align}
	where the last equality holds because $\phi^*$ is independent of $\overline x$. Hence, we have shown that $u_\eps\xrightharpoonup{H^1}u$, with $u(t) = \int_\Om \tilde u(\overline{x},t)\,d\overline{x}$, which concludes the proof of (i).
	
	To see (ii), first use the compact embedding $H^1(\Om)\hookrightarrow L^2(\Om)$ to see that $\bigl\|\eps^{\f{N-1}{2}}\tilde u_\eps-\tilde u\bigr\|_{L^2(\Om)}\to 0$, for a subsequence, and note that $\|\bar\nabla\tilde u_\eps\|_{L^2(\Om)}\to 0$ by assumption. It follows that $\bar\nabla\tilde u=0$, that is $\tilde u(x)=c\cdot u(x_N)$. A simple calculation shows $c=|\Om_0|^{-1}$. Reversing the scaling, this proves (ii).
\end{proof}
In the same way as above one can prove the existence of weakly convergent subsequences in $L^2(\Om_\eps)$.
\begin{prop}
		Let $f_\eps\in L^2(\Om_\eps)$ and $\|f_\eps\|_{L^2(\Om_\eps)}$ uniformly bounded. Then there exists a subsequence $f_{\eps'}$ with $f_{\eps'}\xrightharpoonup{L^2} f$ for some $f\in L^2((0,1))$ as $\eps'\to 0$.
\end{prop}
\begin{proof}
	$L^2$-boundedness in the scaled domain $\Om$ yields weak convergence of $\eps'^{\f{N-1}{2}}f_{\eps'}$ in $L^2(\Om_\eps)$. Scaling back as in the proof of Proposition \ref{weakconv} yields the assertion.
\end{proof}
\section{Proof of Theorem \ref{mainth}}\label{sec:proof_of_theorem}
\subsection{Auxiliary results}
In the following, our discussion will be along the lines of the classical proof from \cite{CM} with the necessary modifications. We define an auxiliary function $w_\eps$ as follows. Let $P_i^\eps$ denote a cube of edge length $2\delta_\eps$ centered at $i\in L_\eps$ and let $w_\eps$ be the solution to 
\begin{equation}\label{wepsilon}
	\begin{cases}
		\phantom{\Delta}w_\eps = 0 & \text{ in } B_{r_\eps}(i),\\
		\Delta w_\eps = 0 &\text{ in } B_{\delta_\eps}(i)\setminus B_{r_\eps}(i),\\
		\phantom{\Delta}w_\eps = 1 &\text{ in }P_i^\eps\setminus B_{\delta_\eps}(i),\\
		\phantom{\Delta}w_\eps &\text{continuous,}
	\end{cases}
\end{equation}
Requiring that $w_\eps\equiv 1$ outside the union of all $P_i^\eps$ we obtain a function $w_\eps\in W^{1,\infty}(\R^N)$ for every $\eps>0$. In fact, exploiting radial symmetry, one can derive the explicit expression
\begin{equation*}
	w_\eps(r)  = \f{r^{2-N}-r_\eps^{2-N}}{\delta_\eps^{2-N}-r_\eps^{2-N}}
\end{equation*}
in polar coordinates (cf. \cite[eq. (2.2)]{CM}). Note that in particular $w_\eps\equiv 1$ in the small cubes $C_j^\eps$ of edge length $\f{2(\sqrt N-1)}{\sqrt N}\delta_\eps$ centered at the corners of the $P_i^\eps$ (cf. \cite[Fig. 2]{CM}).
\begin{lemma}\label{chi}
	Denote $C_\eps:=\bigcup_{j\in L_\eps}C^\eps_j$. The characteristic function $\chi_{C_\eps}$ converges to a constant $\alpha$ weakly$^\star$ in $L^\infty$ in the sense that for all $\varphi\in L^1((0,1))$ and $\varphi_\eps\in L^1(\Om_\eps)$ such that $|\eps\Om_0|^{-1}\|\varphi_\eps-\varphi^*\|_{L^1(\Om_\eps)}\to 0$ as $\eps\to 0$, one has
	$$|\eps\Om_0|^{-1}\langle\chi_{C_\eps},\varphi_\eps\rangle_{L^\infty,L^1}\to \alpha\int_0^1\varphi(x)\,dx$$ 
	(recall the convention $\varphi^*(\overline x,x_N)=\varphi(x_N)$).
\end{lemma}
\begin{proof}
	We use the shorthand $\chi_\eps:=\chi_{C_\eps}$. We first prove the statement for smooth $\varphi$. The general statement will then follow by a density argument. To this end, let $\varphi\in C^\infty((0,1))$ and assume $|\eps\Om_0|^{-1}\|\varphi_\eps-\varphi^*\|_{L^1(\Om_\eps)}\to 0$.  Then
	\begin{align*}
		|\eps\Om_0|^{-1}\int_{\Om_\eps} \chi_\eps\varphi_\eps \,dx &= |\eps\Om_0|^{-1}\int_{\Om_\eps}\chi_\eps \varphi^*\,dx + |\eps\Om_0|^{-1}\int_{\Om_\eps}\chi_\eps(\varphi_\eps-\varphi^*)\,dx \\
		&=: |\eps\Om_0|^{-1}\int_{\Om_\eps}\chi_\eps \varphi^*\,dx + I_\eps.
	\end{align*}
	We have
	\begin{align*}
		|I_\eps| &\leq \|\chi_\eps\|_\infty \cdot |\eps\Om_0|^{-1}\|\varphi_\eps-\varphi^*\|_{L^1(\Om_\eps)}\\
		&\to 0,
	\end{align*}
	by assumption on $\varphi_\eps$. Denote by $x_j^\eps$ the centres of the cubes $C_j^\eps$ and consider the remaining term
	\begin{align*}
		|\eps\Om_0|^{-1}\int_{\Om_\eps}\chi_\eps \varphi^*\,dx &= |\eps\Om_0|^{-1}\sum_j\int_{C_j^\eps} \varphi^*(x_j^\eps)\,dx + |\eps\Om_0|^{-1}\sum_j\int_{C_j^\eps} (\varphi^* - \varphi^*(x_j^\eps))\,dx \\
		&=: |\eps\Om_0|^{-1}\sum_j|C_j^\eps| \varphi^*(x_j^\eps) + \sum_j I^\eps_j.
	\end{align*}
	The total volume of $C_\eps$ is asymptotically $|C_\eps|=\sum_j C_j^\eps\sim |\Om_0| \underbrace{\tfrac{1}{\delta_\eps}\left(\tfrac{\eps}{\delta_\eps}\right)^{N-1}}_{\text{\tiny{number of cubes}}} \underbrace{\phantom{\Big(}\delta_\eps^N\phantom{\Big)}}_{\text{\tiny volume}} = |\eps\Om_0| $. Thus
	\begin{align*}
		\sum_j|I_j^\eps| &\leq |\eps\Om_0|^{-1}\sum_j|C_j^\eps| \|\varphi^* - \varphi^*(x_j^\eps)\|_{L^\infty(C_j^\eps)} \\
		&\leq C \sup_j\|\varphi^* - \varphi^*(x_j^\eps)\|_{L^\infty(C_j^\eps)}\\
		&\to 0\qquad\qquad (\eps\to 0),
	\end{align*}
	where the last statement follows from the smoothness of $\varphi$. Putting the pieces back together we have
	\begin{align*}
		|\eps\Om_0|^{-1}\int_{\Om_\eps} \chi_\eps\varphi_\eps\,dx &= |\eps\Om_0|^{-1}\sum_j |C_j^\eps|\varphi^*(x_j^\eps) + o(1)
	\end{align*}
	Note that the volumes $|C_j^\eps|\sim \delta_\eps^N$ do not depend on $j$ and so
	\begin{align*}
		|\eps\Om_0|^{-1}\int_{\Om_\eps} \chi_\eps\varphi_\eps\,dx &= \alpha'\,\eps^{-N+1}\delta_\eps^N\sum_j \varphi^*(x_j^\eps) + o(1)
	\end{align*}
	for some constant $\alpha'$. Next we use the fact that all $x_j^\eps$ lie in planes $\{x_n=\text{const}\}$ and that $\varphi^*$ is constant in $\bar x$. Thus all terms $\varphi^*(x_j^\eps)$ in the above sum with $(x_j^\eps)_N=(x_k^\eps)_N$ are equal and lead to a factor $\big(\f{\eps}{\delta_\eps}\big)^{N-1}$. Denoting $t_1^\eps,\dots,t_n^\eps$ the projection of $x_j^\eps$ onto the $N$-th coordinate we obtain
	\begin{align*}
		|\eps\Om_0|^{-1}\int_{\Om_\eps} \chi_\eps\varphi_\eps\,dx &= \alpha\,\eps^{-N+1}\delta_\eps^N \left(\f{\eps}{\delta_\eps}\right)^{N-1} \sum_{m=1}^n\varphi(t_m^\eps)  + o(1)\\
		&=\alpha \sum_{m=1}^n \delta_\eps \varphi(t_m^\eps)  + o(1)\\
		&\to \alpha \int_0^1\varphi(t)\,dt
	\end{align*}
	for some constant $\alpha$. The last statement holds because $\varphi$ is Riemann integrable.
	
	Finally we prove the statement for all $\varphi\in L^1((0,1))$. This follows by a standard density argument, though some care is required to deal with the technical difficulties posed by the varying function spaces.
	Let $\varphi\in L^1((0,1))$ be arbitrary and let $\varphi_\eps \in L^1(\Om_\eps)$ such that $|\eps\Om_0|^{-1}\|\varphi_\eps - \varphi^*\|_{L^1(\Om_\eps)} \to 0$.
	Next, let $\delta>0$ and use density of $C^\infty((0,1))$  in $L^1((0,1))$ to choose $\eta\in C^\infty((0,1))$ with $\| \varphi_ - \eta \|_{L^1((0,1))} < \delta$ and let $\eta_\eps\in L^1(\Om_\eps)$ be such that $|\eps\Om_0|^{-1}\|\eta_\eps - \eta^*\|_{L^1(\Om_\eps)} \to 0$. We first note that $\varphi_\eps$ and $\eta_\eps$ are necessarily close in the limit:
	\begin{align}
		\limsup_{\eps\to 0} |\eps\Om_0|^{-1}\|\varphi_\eps	- \eta_\eps\|_{L^1(\Om_\eps)} 
		&\leq \limsup_{\eps\to 0} |\eps\Om_0|^{-1}\big( \|\varphi_\eps - \varphi^*\|_{L^1(\Om_\eps)} + \|\varphi^* - \eta^*\|_{L^1(\Om_\eps)} + \|\eta^* - \eta_\eps\|_{L^1(\Om_\eps)} \big) \nonumber \\
		&\leq \limsup_{\eps\to 0} |\eps\Om_0|^{-1}\|\varphi^* - \eta^*\|_{L^1(\Om_\eps)} \nonumber  \\
		&= \|\varphi-\eta\|_{L^1((0,1))}, \nonumber \\
		&< \delta \label{eq:density1}
	\end{align}
	where the second line follows from the assumptions on $\eta_\eps$ and $\varphi_\eps$ and the third line follows from the definition of $\varphi^*$ and $\eta^*$.
	Next, we estimate
	\begin{align*}
		\left| |\eps\Om_0|^{-1}\langle\chi_\eps,\varphi_\eps\rangle - \alpha\int_0^1 \varphi(t)\,dt \right|
		&\leq |\eps\Om_0|^{-1} |\langle \chi_\eps, \varphi_\eps - \eta_\eps\rangle| \\
		&\qquad + \left| |\eps\Om_0|^{-1}\langle\chi_\eps,\eta_\eps\rangle - \alpha\int_0^1 \eta(t)\,dt \right|\\
		&\qquad + |\alpha|\int_0^1 |\eta(t)-\varphi(t)|\,dt\\
		&\leq \|\chi_\eps\|_\infty |\eps\Om_0|^{-1}\|\varphi_\eps - \eta_\eps\|_{L^1(\Om_\eps)} \\
		&\qquad + \left| |\eps\Om_0|^{-1}\langle\chi_\eps,\eta_\eps\rangle - \alpha\int_0^1 \eta(t)\,dt \right|\\
		&\qquad + |\alpha|\delta.
	\end{align*}
	Finally, using \eqref{eq:density1}, together with the facts that $\|\chi_\eps\|_\infty\leq 1$ and $|\eps\Om_0|^{-1}\langle\chi_\eps,\eta_\eps\rangle \to \alpha\int_0^1 \eta(t)\,dt $, we conclude that
	\begin{align*}
		\limsup_{\eps\to 0} \left| |\eps\Om_0|^{-1}\langle\chi_\eps,\varphi_\eps\rangle - \alpha\int_0^1 \varphi(t)\,dt \right| 
		\leq (1 + |\alpha|)\delta.
	\end{align*}
	Since $\delta>0$ was arbitrary, it follows that
	\begin{align*}
		\limsup_{\eps\to 0} \left| |\eps\Om_0|^{-1}\langle\chi_\eps,\varphi_\eps\rangle - \alpha\int_0^1 \varphi(t)\,dt \right| =0.
	\end{align*}
\end{proof}
\begin{lemma}\label{w=1}
	For the function $|\eps\Om_0|^{-\f 1 2}w_\eps$, with $w_\eps$ defined in \eqref{wepsilon}, one has	$|\eps\Om_0|^{-\f 1 2}w_\eps\xrightharpoonup{H^1}1$. 
\end{lemma}
\begin{proof}
	It follows by a trivial modification of the argument in \cite{CM} that $|\eps\Om_0|^{-\f 1 2}w_\eps$ satisfies the bound \eqref{H1bdd} and even the stronger condition (ii) in Proposition \ref{weakconv}. Thus, by Proposition \ref{weakconv} there exists a subsequence $|\eps\Om_0|^{-\f 1 2}w_\eps\xrightharpoonup{H^1} w$ for some $w\in H^1((0,1))$ and $|\eps\Om_0|^{-\f 1 2}w_\eps\xrightarrow{L^2}w$. It remains to show $w=1$. This will be done by applying Lemma \ref{chi}.
	\begin{itemize}
		\item[\emph{Claim:}] If $\phi_\eps\xrightarrow{L^2} \phi$ then $|\eps\Om_0|^{-1}\bigl\|w_\eps|\eps\Om_0|^{\f 1 2}\phi_\eps-w^*\phi^*\bigr\|_{L^1(\Om_\eps)}\to 0$.
		\item[\emph{Proof of claim:}] By the triangle inequality we have
		\begin{align*}
			\hspace{-8mm}|\eps\Om_0|^{-1}\left\|w_\eps|\eps\Om_0|^{\f 1 2}\phi_\eps - w^*\phi^*\right\|_{L^1(\Om_\eps)} &\leq |\eps\Om_0|^{-1}\left\|w_\eps|\eps\Om_0|^{\f 1 2}\phi_\eps - w_\eps\phi^*\right\|_{L^1(\Om_\eps)} \\
			&\qquad\qquad + |\eps\Om_0|^{-1}\bigl\|w_\eps\phi^* - w^*\phi^*\bigr\|_{L^1(\Om_\eps)}\\
			&\leq |\eps\Om_0|^{-1}\|w_\eps\|_{L^2(\Om_\eps)} \left\||\eps\Om_0|^{\f 1 2}\phi_\eps - \phi^*\right\|_{L^2(\Om_\eps)} \\
			&\qquad\qquad + |\eps\Om_0|^{-1}\|\phi^*\|_{L^2(\Om_\eps)}\bigl\|w_\eps - w^*\bigr\|_{L^2(\Om_\eps)}\\
			&= \left( |\eps\Om_0|^{-\f12}\|w_\eps\|_{L^2(\Om_\eps)} \right) \left(\left\|\phi_\eps - U_\eps\phi\right\|_{L^2(\Om_\eps)}\right)\\
			&\qquad + \left( |\eps\Om_0|^{-\f12}\|\phi^*\|_{L^2(\Om_\eps)} \right)\left(\left\||\eps\Om_0|^{-\f 1 2}w_\eps - U_\eps w \right\|_{L^2(\Om_\eps)}\right)\\
			&\to 0
		\end{align*}
	\end{itemize}
To prove $w=1$, note that $w_\eps\chi_{C_\eps} = \chi_{C_\eps}$. Hence, for $\phi_\eps\xrightarrow{L^2}\phi$ Lemma \ref{chi} (with $\varphi_\eps=w_\eps|\eps\Om_0|^{\f 1 2}\phi_\eps$) gives
	\begin{align*}
		|\eps\Om_0|^{-\f 1 2}\int_{\Om_\eps} w_\eps \chi_{C_\eps}\phi_\eps\,dx 
		&= 
		|\eps\Om_0|^{-1}\int_{\Om_\eps} \underbrace{w_\eps|\eps\Om_0|^{\f 1 2}\phi_\eps}_{\text{\tiny str. in }L^1}\chi_{C_\eps}\,dx\\
		&\to \; \alpha \int_0^1 w\phi\,dx.
	\end{align*}
	On the other hand, also by Lemma \ref{chi}
	\begin{align*}
		|\eps\Om_0|^{-\f12}\int_{\Om_\eps} \chi_{C_\eps}\phi_\eps\,dx &= |\eps\Om_0|^{-1}\int_{\Om_\eps} \chi_{C_\eps}|\eps\Om_0|^{\f12}\phi_\eps\,dx\\
		&\to \; \alpha \int_0^1 \phi\,dx.
	\end{align*}
	Since $\phi\in L^2((0,1))$ was arbitrary, we conclude $w=1$. 
\end{proof}
From \hyperref[w=1]{Lemma} \ref{w=1} we conclude that $|\eps\Om_0|^{-\f 1 2}\nabla w_\eps\xrightharpoonup{L^2} 0$ (note that this is the full gradient and not merely $\bar\nabla$), i.e. we have
\begin{align}\label{gradwto0}
	\int_{\Om_\eps}|\eps\Om_0|^{-\f12}\nabla w_\eps\cdot\boldsymbol\psi_\eps\,dx\to 0 
\end{align}
whenever $\|\boldsymbol\psi_\eps-U_\eps\boldsymbol\psi\|_{L^2(\Om_\eps)^N}\to 0$ for some $\boldsymbol\psi\in L^2((0,1))^N$.
\subsection{Convergence of solutions}\label{sec:sonvergence_of_solutions}
\begin{lemma}
	Let $u_\eps$ be a weak solution of \eqref{pde} with right hand side $f_\eps\xrightarrow{L^2}f$. Then the a priori bound
	\begin{align}\label{apriori}
		 \|u_\eps\|_{L^2(\Om_\eps)}^2 + \|\nabla u_\eps\|_{L^2(\Om_\eps)}^2 \leq C \|f\|_{L^2((0,1))}^2
	\end{align}
	holds.
\end{lemma}
\begin{proof}
	The weak formulation of \eqref{pde} yields for arbitrary $\delta>0$
	\begin{align*}
		\int_{\Om_\eps} |\nabla u_\eps|^2\,dx + z\int_{\Om_\eps}|u_\eps|^2\,dx &= \int_{\Om_\eps} f_\eps u_\eps\,dx \nonumber\\
		&\leq \f{\delta}{2} \|u_\eps\|_{L^2(\Om_\eps)}^2 + (2\delta)^{-1}\|f_\eps\|_{L^2(\Om_\eps)}^2
	\end{align*}
	Choosing e.g. $\delta:=z$ this yields
	\begin{align}\label{eq:apriori_proof}
		\|\nabla u_\eps\|_{L^2(\Om_\eps)}^2 + \f{z}{2}\|u_\eps\|_{L^2(\Om_\eps)}^2 &\leq  (2z)^{-1}\|f_\eps\|_{L^2(\Om_\eps)}^2
	\end{align}
		Next, {without loss of generality,} choose $\eps$ small enough such that $\big|\|f_\eps\|_{L^2(\Om_\eps)}^2 - \|f\|_{L^2((0,1))}^2\big|<\|f\|_{L^2((0,1))}^2$. We obtain from \eqref{eq:apriori_proof}
	\begin{align*}
		\|\nabla u_\eps\|_{L^2(\Om_\eps)}^2 + \f{z}{2}\|u_\eps\|_{L^2(\Om_\eps)}^2 &\leq  \left((2z)^{-1}+1 \right)\|f\|_{L^2((0,1))}^2
	\end{align*}
	and hence
	\begin{align*}
		\|\nabla u_\eps\|_{L^2(\Om_\eps)}^2 + \|u_\eps\|_{L^2(\Om_\eps)}^2 &\leq  \f{(2z)^{-1}+1}{\min\{1,\nicefrac{z}{2}\}} \|f\|_{L^2((0,1))}^2.
	\end{align*}
\end{proof}
Note that this a priori bound actually proves that case (ii) of Lemma \ref{weakconv} is satisfied by the solutions $u_\eps$, since $\|\bar\nabla u_\eps\|_{L^2(\Om_\eps)}$ is uniformly bounded. Thus there exists $u\in H^1((0,1))$ such that $u_\eps\xrightharpoonup{H^1}u$ and $u_\eps\xrightarrow{L^2}u$. We will show that $u$ satisfies the weak version of \eqref{limeq}. Let $\phi\in H^1((0,1))$ and consider the weak formulation of \eqref{pde} with test function $w_\eps\cdot U_\eps\phi$:
\begin{align}\label{weakform}
	\int_{\Om_\eps} \overline{\nabla u}_\eps\cdot\nabla(w_\eps U_\eps\phi)\,dx + z\int_{\Om_\eps}\overline{u}_\eps w_\eps U_\eps\phi\,dx &= \int_{\Om_\eps}\overline{f}_\eps w_\eps U_\eps\phi\,dx\nonumber\\
	\Leftrightarrow\; \int_{\Om_\eps} (U_\eps\phi)\,\overline{\nabla u}_\eps\cdot\nabla w_\eps \,dx + \int_{\Om_\eps} w_\eps\overline{\nabla u}_\eps\cdot\nabla (U_\eps\phi) \,dx
	+z\int_{\Om_\eps}\overline{u}_\eps w_\eps U_\eps\phi\,dx &= \int_{\Om_\eps}\overline{f}_\eps w_\eps U_\eps\phi\,dx
\end{align}
We will consider the convergence {all four} terms separately. 
\paragraph{Right-hand side:}
Since $\phi\in H^1((0,1))$ we have $\|\phi\|_{L^\infty}<C\|\phi\|_{H^1((0,1))}$ uniformly in $\eps$, by Morrey's inequality. Thus $w_\eps U_\eps\phi$ converges strongly in $L^2$ to $\phi$. Indeed, we have 
	\begin{align*}
		\|w_\eps U_\eps\phi-U_\eps\phi\|_{L^2(\Om_\eps)} &\leq \|U_\eps\phi\|_{\infty}\|w_\eps -1\|_{L^2(\Om_\eps)}\\
		&= \|\phi\|_\infty \left\||\eps\Om_0|^{-\f12}w_\eps - U_\eps(1)\right\|_{L^2(\Om_\eps)}\\
		&\to 0
	\end{align*}
Since $f_\eps\xrightharpoonup{L^2} f$ we can conclude
\begin{align*}
	\int_{\Om_\eps}\overline f_\eps w_\eps U_\eps\phi\,dx\,\to\, \int_0^1 \overline f\phi\,dx
\end{align*}
\paragraph{Third term on the l.h.s.:}
By the same reasoning as above, one has $u_\eps\to u$ and $w_\eps U_\eps\phi\to \phi$ strongly in $L^2$ and thus
\begin{align*}
	z\int_{\Om_\eps}\overline u_\eps w_\eps U_\eps\phi\,dx\,\to\, z\int_0^1 \overline u\phi\,dx
\end{align*}
\paragraph{Second term on the l.h.s.:}
By the same reasoning as above, $w_\eps \nabla (U_\eps\phi)=w_\eps U_\eps\phi'$ converges strongly in $L^2$ to $\phi'$. Since $\nabla u_\eps$ converges weakly  in $L^2$, the whole integral converges to $\int_0^1\overline u'\phi'\,dt$.
\paragraph{First term on the l.h.s.:}
First, we rewrite the term
\begin{align}\label{str}
	\int_{\Om_\eps}(U_\eps\phi)\,\overline {\nabla u}_\eps\cdot\nabla w_\eps dx &= \langle-\Delta w_\eps , u_\eps U_\eps\phi\rangle_{H^{-1},H^1_0} - \int_{\Om_\eps} \overline u_\eps\nabla w_\eps\cdot\nabla (U_\eps\phi)\,dx
\end{align}
The second term on the right hand side of \eqref{str} converges to 0 by \eqref{gradwto0}. Indeed, since $u$ and $\nabla U_\eps\phi$ are uniformly bounded in $L^\infty$, by Morrey's inequality, we have $\overline u_\eps\nabla U_\eps\phi\xrightarrow{L^2}\overline u \phi'$.
The last remaining term is treated in the following
\begin{lemma}\label{lemma:laplace_w}
	One has
	\begin{align*}
		\langle-\Delta w_\eps , u_\eps U_\eps\phi\rangle_{H^{-1},H^1_0}\;\to\; \mu\int_0^1 \overline u\phi\,dt,
	\end{align*}
	where $\mu$ was defined Theorem \ref{mainth}. 
\end{lemma}
\begin{proof}
	The proof is only a small variation of that of \cite[Lemma 2.3]{CM}. We give it here nevertheless for the sake of self-containedness. First, note that by partial integration and boundary conditions, we have
	$$
		\langle-\Delta w_\eps , u_\eps\phi_\eps\rangle = \f{N-2}{1-\delta_\eps^2} \sum_{i\in L_\eps} \langle S_i^\eps,u_\eps U_\eps\phi\rangle,
	$$
	where $S_i^\eps$ is the Dirac measure on $\del B_{\delta_\eps}(i)$: $\langle S_i^\eps , \varphi\rangle = \int_{\del B_{\delta_\eps}(i)}\varphi\,dS$. Moreover, let us define the function $q_\eps$ as the unique solution of the Neumann problem
	\begin{align*}
		\begin{cases}
			-\Delta q_\eps = N, & \text{in }B_{\delta_\eps}(i)\\
			\;\;\, \del_\nu  q_\eps = \eps & \text{on }\del B_{\delta_\eps}(i)
		\end{cases}
	\end{align*}
	satisfying $q_\eps=0$ on $\del B_{\delta_\eps}(i)$. Extending $q_\eps$ by zero to all of $\Om_\eps$ we can easily see that $q_\eps \to 0$ in $W^{1,\infty}(\R^N)$. Consequently:
	\begin{align*}
		\langle -\Delta q_\eps,\varphi_\eps \rangle &=  \int_{\Om_\eps} \nabla q_\eps\cdot\nabla\varphi_\eps\,dx\\
		&\leq \|\nabla q_\eps\|_{\infty}\cdot \|\varphi_\eps\|_{L^1(\Om_\eps)}\\
		&\to 0,
	\end{align*}
	for every sequence with $\|\varphi_\eps\|_{L^1(\Om_\eps)}$ bounded. On the other hand, one has $-\Delta q_\eps = N\chi^\eps_{\cup_i B_{\delta_\eps}(i)} - \sum_{i\in L_\eps}\delta_\eps S_i^\eps$. Thus, we can take the limit in the following equation
	\begin{align*}
		\langle -\Delta q_\eps,\varphi_\eps \rangle &= \int\limits_{\cup_i B_{\delta_\eps}(i)}\varphi_\eps\,dx + \sum_{i\in L_\eps}\,\delta_\eps\!\int\limits_{\del B_{\delta_\eps}(i)} \varphi_\eps\,dS.
	\end{align*} 
	The first term on the right hand side converges to $\mu\int_0^1 u\phi\,dt$ as can be seen by the same argument as in the proof of Lemma \ref{chi}. We obtain the equality
	\begin{align*}
		\lim_{\eps\to 0} \sum_{i\in L_\eps}\,\delta_\eps\!\int\limits_{\del B_{\delta_\eps}(i)} \varphi_\eps\,dS = \mu\int_0^1 \varphi\,dt.
	\end{align*}
	The assertion now follows by choosing $\varphi_\eps = u_\eps U_\eps\phi$ in the above equation (note that $\|u_\eps U_\eps\phi\|_{L^1(\Om_\eps)}$ is uniformly bounded).
\end{proof}
This settles the convergence of the last remaining term in \eqref{weakform} and leads to the limit problem
\begin{align}\label{limitprob}
	\int_0^1 \overline u'\phi'\,dt + (\mu+z)\int_0^1 \overline u\phi\,dt = \int_0^1 \overline f\phi\,dt,
\end{align}
which is nothing but the weak formulation of \eqref{limeq}. Since it has already been shown that $u_\eps$ satisfies hypothesis (ii) of Proposition \ref{weakconv} and thus converges strongly in $L^2$, the proof of Theorem \ref{mainth} is completed.
\qed
\begin{remark}
	We note that our assumption on the spherical shape of the holes was made for the sake of definiteness, however, our results easily generalise to more general geometries as detailed in \cite[Th. 2.7]{CM}. Moreover, our results are also valid  for more general elliptic operators $\operatorname{div}(A\nabla)$ with continuous coefficients $A$ (cf. \cite[Ex. 2.16]{CM}).
\end{remark}

\section{Norm-Resolvent Convergence}\label{NormandSpec}

In this section we will take a more operator-theoretic point of view and prove operator norm convergence for the resolvent. To this end, let us first introduce some notation. We define the following operators in $L^2$.
\begin{equation}\label{Adef}
\begin{aligned}
	A_\eps & :=-\Delta,  &\;\mathcal D(A_\eps)&=\{u\in\mathcal H_\eps^0\cap H^2(\Oeps) : \del_\nu u|_{\del\Om_\eps}=0\}\\
	A  & := -\f{d^2}{dt^2}+\mu, &\;\mathcal D(A)&=\{ u\in H^2((0,1)) :  u'(0)= u'(1)=0\},
\end{aligned}
\end{equation}
{where $\mathcal D(\cdot)$ denotes the domain of the relevant operator.} Furthermore, we define the two identification operators between the domains
\begin{equation}\label{UVdef}
\begin{aligned}
	&\U:L^2((0,1)) \to L^2(\Oeps); \;&\; (\U g)(x) &= |\eps\Om_0|^{-\f12}g(x_N)\\
	&\Vt:L^2(\Oeps) \to L^2((0,1)); \;&\; (\Vt f)(t) &= |\eps\Om_0|^{-\f12}\int_{\eps\Om_0}\widetilde f(\bar x,t)\,d\bar x,
\end{aligned}
\end{equation}
where $\widetilde f$ denotes extension { of $f$} by 0 into the holes. Note that $\|\U\|_{\mathcal L(L^2((0,1)),L^2(\Oeps))},\|\Vt\|_{\mathcal L(L^2(\Oeps),L^2((0,1)))}$ are uniformly bounded in $\eps$.

Now, let us go back to \eqref{weakform}, and observe that the right-hand side will still converge if $f_\eps$ is only \emph{weakly} convergent in $L^2$. We deduce the following
\begin{lemma}\label{uniformU}
	Let $(g_\eps)\subset L^2((0,1))$ and assume that $g_\eps\rightharpoonup g$ weakly in $L^2((0,1))$. Then for any $z>0$ one has
	\begin{align*}
		\|(A_\epsilon+z)^{-1}\U g_\eps - \U(A+z)^{-1}g\|_{L^2(\Oeps)}\to 0
	\end{align*}
	in $L^2((0,1))$.
\end{lemma}
\begin{proof}
	By the above comment, it is enough to show that $\U g_\eps\xrightharpoonup{L^2}g$ in the sense of Definition \ref{weakdef}. To this end, let $\phi_\eps\in L^2(\Oeps)$ and assume $\phi_\eps\xrightarrow{L^2}\phi$ for some $\phi\in L^2((0,1))$. We have
	\begin{align*}
		\langle\U g_\eps,\phi_\eps\rangle_{L^2(\Oeps)} &= 
	\langle\U g_\eps,\U\phi\rangle_{L^2(\Oeps)} + \langle\U g_\eps,\phi_\eps - \U\phi\rangle_{L^2(\Oeps)}\\
	&= \langle\U g_\eps,\U\phi\rangle_{L^2(\Om_\eps)} + \langle\U g_\eps,\U\phi\rangle_{L^2(T_\eps)} + \langle\U g_\eps,\phi_\eps - \U\phi\rangle_{L^2(\Oeps)}\\
	&= \langle g_\eps,\phi\rangle_{L^2((0,1))} + \langle\U g_\eps,\U\phi\rangle_{L^2(T_\eps)} + \langle\U g_\eps,\phi_\eps - \U\phi\rangle_{L^2(\Oeps)}
	\end{align*}
	The last term goes to 0 since $\phi_\eps\xrightarrow{L^2}\phi$, wheres the second term on the right hand side converges to 0 because $|\eps^{-1}T_\eps|\to 0$. Finally, the first term on the right-hand side converges to $\langle g,\phi\rangle_{L^2((0,1))}$ by assumption, which concludes the proof.
\end{proof}
Lemma \ref{uniformU} shows that using $\U$ as an identification operator, the convergence of solutions of \eqref{pde} is uniform in the right-hand side. We will now prove a similar statement for $\Vt$. 
\begin{lemma}\label{uniformV}
	Let $f_\eps\in L^2(\Oeps)$ be a sequence with $f_\eps\xrightharpoonup{L^2}f$ and $u_\eps$ be the sequence of solutions to \eqref{pde}. Then one has 
	\begin{align*}
		\Vt u_\eps\rightharpoonup u\quad \text{ in } H^1((0,1)),
	\end{align*}
	where $u$ solves the limit problem \eqref{limitprob}.
\end{lemma}
\begin{proof}
	First, note that $\|\Vt u_\eps\|_{H^1((0,1))}$ is uniformly bounded in $\eps$. Indeed, we can compute
	\begin{align*}
		\|\Vt u_\eps\|_{H^1((0,1))}^2 
		&\stackrel{\phantom{\text{Jensen}}}{=} 
		 \int_0^1\left| |\eps\Om_0|^{-\f12}\int_{\eps\Om_0}u_\eps(\bar x,t) \,d\bar x\right|^2\,dt + \int_0^1\left||\eps\Om_0|^{-\f12} \int_{\eps\Om_0}\del_N u_\eps(\bar x,t) \,d\bar x\right|^2\,dt\\
		&\stackrel{\text{Jensen}}{\leq} 
		 \int_0^1 \int_{\eps\Om_0}\left|u_\eps(\bar x,t) \right|^2\,d\bar x dt +  \int_0^1 \int_{\eps\Om_0}\left|\del_N u_\eps(\bar x,t) \right|^2\,d\bar x dt\\
		&\stackrel{\phantom{\text{Jensen}}}{\leq}
		\|u_\eps\|_{L^2(\Oeps)}^2 + \|\nabla u_\eps\|_{L^2(\Oeps)}^2 \\
		&\stackrel{\phantom{\text{Jensen}}}{\leq} C\|f_\eps\|^2_{L^2(\Oeps)}
	\end{align*}
	by the a priori bound \eqref{apriori}. The right hand side remains bounded as $\eps\to 0$ since $(f_\eps)$ converges weakly. Hence there exists a $H^1$-weakly convergent subsequence (again denoted by $\Vt u_\eps$) with $\Vt u_\eps\rightharpoonup v$ for some $v\in H^1((0,1))$.
	By the Rellich-Kondrachov theorem one has $\Vt u_\eps\rightarrow v$ strongly in $L^2((0,1))$. It remains to show that $v=u$. This will be done in two steps. Step 1: Because $f_\eps\rightharpoonup f$, every term in the weak formulation \eqref{weakform} converges, that is, $u_\eps\xrightharpoonup{H^1} u$ (and thus strongly in $L^2$) in the sense of Definition \ref{weakdef}, where $u$ solves the limit problem \eqref{limitprob}. Step 2: compute
	\begin{align*}
		\|\Vt u_\eps - u\|_{L^2((0,1))}^2 
		&=
		\int_0^1 \left| |\eps\Om_0|^{-\f12}\int_{\eps\Om_0} u_\eps(\bar x,t)\,d\bar x - |\eps\Om_0|^{-\f12}u(t) \right|^2\,dt\\
		&=
		 \int_0^1 \left| |\eps\Om_0|^{-\f12}\int_{\eps\Om_0} \big(u_\eps(\bar x,t) - |\eps\Om_0|^{-\f12}u(t)\big) \,d\bar x\right|^2\,dt\\
		&\leq 
		\int_0^1  \int_{\eps\Om_0} \left|u_\eps(\bar x,t) - |\eps\Om_0|^{-\f12}u(t)\right|^2 \,d\bar x\,dt\\
		&=
		C\left\|u_\eps - \U u\right\|^2_{L^2(\Om_\eps)}\\
		&\to 0,
	\end{align*}
	where the third line follows from Jensen's inequality, and thus $\Vt u_\eps\to u$ in $L^2((0,1))$ which implies $v=u$ and concludes the proof.
\end{proof}
We are now able to state the main result of this section.
\begin{theorem}\label{normconv}
	Let $A_\eps,A$ and $\U,\Vt$ be defined as in \eqref{Adef} and \eqref{UVdef}. Then one has
	\begin{align}
		\left\| (A_\eps+z)^{-1}\U - \U (A+z)^{-1} \right\|_{\mathcal L(L^2((0,1)),L^2(\Oeps))} &\to 0	\label{Uconv}\\
		\nonumber\\[-2mm]
		\left\| \Vt(A_\eps+z)^{-1} -  (A+z)^{-1}\Vt \right\|_{\mathcal L(L^2(\Oeps),L^2((0,1)))} &\to 0.	\label{Vconv}
	\end{align}
\end{theorem}
\begin{proof}
	We first prove \eqref{Uconv}. Let $(g_\eps)$ be any bounded sequence in $L^2((0,1))$. Then there exists a weakly convergent subsequence $g_{\eps'}\rightharpoonup g$ for some $g\in L^2((0,1))$. Now compute
	\begin{align*}
		\left\| (A_{\eps'}+z)^{-1}\mathcal U_{\eps'} g_{\eps'} - \mathcal U_{\eps'} (A+z)^{-1}g_{\eps'} \right\|_{L^2(\Om_{\eps'}^{\mathrm p})} &\leq 
		\left\| (A_{\eps'}+z)^{-1}\mathcal U_{\eps'}g_{\eps'} - \mathcal U_{\eps'} (A+z)^{-1}g \right\|_{L^2(\Om_{\eps'}^{\mathrm p})} \\
		&\qquad\qquad\qquad + \left\| \mathcal U_{\eps'} (A+z)^{-1}(g - g_{\eps'}) \right\|_{L^2(\Om_{\eps'}^{\mathrm p})}.
	\end{align*}
	The first term on the right hand side converges to 0 by Lemma \ref{uniformU}. The second term converges to 0 too, because $g_{\eps'}\rightharpoonup g$, $(A+z)^{-1}$ is a compact operator and $\|\U\|_{\mathcal L(L^2((0,1)),L^2(\Oeps))}$ is uniformly bounded. Next, choose $(g_\eps)$ {with $\|g_\eps\|_{L^2((0,1))}\leq 1$} in such a way that 
	$$\sup_{\|h\|_{L^2((0,1))}\leq 1}\left\| \Big((A_\eps+z)^{-1}\U  - \U (A+z)^{-1}\Big)h \right\|_{L^2(\Oeps)}-\eps < \left\| (A_\eps+z)^{-1}\U g_\eps - \U (A+z)^{-1}g_\eps \right\|_{L^2(\Oeps)}.$$
	By the above, the right-hand side of this equation converges to 0 for a suitable subsequence $(\eps')$, so taking the limit $\eps'\to 0$ on both sides yields
	\begin{align*}
		\limsup_{\eps'\to 0}\sup_{\|h\|_{L^2((0,1))}\leq 1}\left\| \Big((A_{\eps'}+z)^{-1}\mathcal U_{\eps'}  - \mathcal U_{\eps'} (A+z)^{-1}\Big)h \right\|_{L^2(\Om_{\eps'}^{\mathrm p})} \leq 0.
	\end{align*}
	Applying this reasoning to every subsequence of $(A_\eps+z)^{-1}\U - \U (A+z)^{-1}$ yields the claim for the whole sequence and concludes the proof of \eqref{Uconv}.
	
	To prove \eqref{Vconv}, let $f_\eps\in L^2(\Oeps)$ be a sequence with $\|f_\eps\|_{L^2(\Oeps)}$ uniformly bounded. Then there exists $f\in L^2((0,1))$ and a weakly convergent subsequence $(f_{\eps'})$ such that $\widetilde f_{\eps'}\xrightharpoonup{L^2} f$ in the sense of Definition \ref{weakdef} (where $\widetilde f_\eps$ denotes extension by 0 from $\Oeps$ to $\Om_\eps$). In particular we have 
	$$\int_{\Om_{\eps'}}\widetilde f_{\eps'} \mathcal U_{\eps'}\phi\,dx = \int_{\Om_{\eps'}^{\mathrm p}} f_{\eps'} \mathcal U_{\eps'}\phi\,dx\to \int_0^1 f\phi, dt $$
	as $\eps'\to 0$. The left hand side  of this equation can be rewritten in terms of $\Vt f_\eps$:
	\begin{align*}
		\int_{\Oeps}f_{\eps}\, \U\phi\,dx &= \int_0^1\int_{\eps\Om_0}|\eps\Om_0|^{-\f12}\widetilde f_\eps(\bar x,t)\,d\bar x\,\phi(t)\,dt\\
		&= \int_0^1 (\Vt f_\eps)\phi\,dt.
	\end{align*}
	Hence we have $\tilde{\mathcal U}_{\eps'} f_{\eps'}\rightharpoonup f$ in $L^2((0,1))$. The rest of the proof is entirely analogous to that of \eqref{Uconv}, using compactness of $(A+z)^{-1}$ and Lemma \ref{uniformV}.
\end{proof}
\section{Spectral Convergence}\label{SpSec}
In this section we will prove Corollary \ref{SpCon}. Let us first note that, since the domains $\Oeps$ and $(0,1)$ are bounded, the domains $\mathcal D(A_\eps), \mathcal D(A)$ are compactly embedded in $L^2$ and hence $A_\eps$ and $A$ have compact resolvent and their spectra are discrete. Let us denote by $(\lambda_k^\eps)$, resp. $(\lambda_k)$, the eigenvalues of $A_\eps+\mathrm{id}$, resp. $A+\mathrm{id}$, labelled in increasing order. We will use a theorem from \cite{IOS} to prove the convergence of spectra.
\begin{theorem}[{\cite[Th. III.1.4]{IOS}}]\label{IOSThm}
	Assume that the following hypotheses are satisfied:
	\begin{itemize}
	\item[\emph{(H1)}] One has $\|\U g\|_{L^2(\Oeps)}\to \|g\|_{L^2((0,1))}$ for all $g\in L^2((0,1))$;
	\item[\emph{(H2)}] The operators $(A_\eps+\mathrm{id})^{-1},(A+\mathrm{id})^{-1}$ are positive, compact, self-adjoint and $\|(A_\eps+\mathrm{id})^{-1}\|_{\mathcal L(L^2(\Oeps))}$ is uniformly bounded in $\eps$;
	\item[\emph{(H3)}] For any $g\in L^2((0,1))$ one has $\| (A_\eps+\mathrm{id})^{-1}\U g - \U(A+\mathrm{id})^{-1} g\|_{L^2(\Oeps)}\to 0$ as $\eps\to 0$
	\item[\emph{(H4)}] For each $f_\eps\in L^2(\Oeps)$ with $\|f_\eps\|_{L^2(\Oeps)}$ uniformly bounded there exists a subsequence $f_{\eps'}$ and some $g\in L^2((0,1))$ such that $\|(A_{\eps'}+\mathrm{id})^{-1}f_{\eps'} - \mathcal U_{\eps'}g\|_{L^2(\Om_{\eps'}^\mathrm{p})}\to 0$ as $\eps'\to 0$.
	\end{itemize}
	Then there exists $C>0$ such that 
	\begin{align}\label{SpEst}
		\left| (\lambda_k^\eps)^{-1} - \lambda_k^{-1}\right| \leq C \sup_{\substack{g\in\operatorname{Eig}(A_0;\lambda_k) \\ \|g\|_{L^2}=1 }}\left\| (A_{\eps'}+\mathrm{id})^{-1}\U g - \U(A+\mathrm{id})^{-1}g \right\|_{\mathcal L(L^2(\Oeps))}
	\end{align}
\end{theorem}
We remark that the constant $C$ in \eqref{SpEst} can be given explicitly in terms of the $\lambda_k$. This more precise version of \eqref{SpEst} can be found in \cite[Eq. (III.1.13)]{IOS}.

We will now show that (H1)-(H4) are satisfied for $A_\eps,A$ and $\U$. First, note that (H2) is obvious from the preceding discussion and the a priori estimate \eqref{apriori}. Furthermore, (H3) follows directly from Theorem \ref{normconv}. (H4) can be seen as follows. If $\|f_\eps\|_{L^2(\Oeps)}\leq C$, there exists a subsequence $f_{\eps'}\xrightarrow{L^2} f$ for some $f\in L^2((0,1))$. Now go back to the weak formulation \eqref{weakform} and note that the right-hand side term $\int_{\Om_{\eps'}}f_{\eps'} w_{\eps'} \mathcal U_{\eps'}\phi\,dx$ only requires \emph{weak} convergence of $f_\eps$ in order to yield the desired limit. This shows (H4) with $g=\big(-\f{d^2}{dt^2}+1+\mu\big)^{-1}f$. Finally, let us prove (H1). We have
\begin{align*}
	\|\U g\|_{L^2(\Oeps)}^2 &= \int_{\Oeps}|\eps\Om_0|^{-1} |g(x_N)|^2\,dx\\
	&= \int_{\Om_\eps}|\eps\Om_0|^{-1} |g(x_N)|^2\,dx + \int_{T_\eps}|\eps\Om_0|^{-1} |g(x_N)|^2\,dx\\
	&= \int_0^1 |g(t)|^2\,dt +\int_{\eps^{-1}T_\eps}|\Om_0|^{-1} |g(x_N)|^2\,dx\\
	&\to \int_0^1 |g(t)|^2\,dt.
\end{align*}
Indeed, one has $|\eps^{-1}T_\eps| \sim \eps^{-N+1}r_\eps^N\f{\eps^{N-1}}{\delta_\eps^N} = \delta_\eps^{\f{2N}{N-2}}\to 0$ as $\eps\to 0$. 

Thus, all hypotheses are satisfied and Theorem \ref{IOSThm} applies. From \eqref{SpEst} we immediately obtain
\begin{align}
	\left| (\lambda_k^\eps)^{-1} - \lambda_k^{-1}\right| \leq C \left\| (A_\eps+z)^{-1}\U - \U (A+z)^{-1} \right\|_{\mathcal L(L^2((0,1)),L^2(\Oeps))}.
\end{align}
Clearly, denoting $a(\eps):= \left\| (A_\eps+z)^{-1}\U - \U (A+z)^{-1} \right\|_{\mathcal L(L^2((0,1)),L^2(\Oeps))}$, this proves Corollary \ref{SpCon}.
\qed
\begin{remark}\label{remark}
	Let us note that all the above results also hold in two dimensions with minor modifications in the definition of the function $w_\eps$ which are detailed in \cite{CM}. We have excluded this case merely to simplify the presentation.
\end{remark}

\section{Graph-like Domains}\label{sec:graph-like_Domains}
In this section we extend our analysis towards domains approximating not merely an interval, but a finite connected graph. That is, the perforated domain consists of ``fattened edges'' of the form $E_\eps:=\eps\Om_0\times(0,\ell)$ which are connected by ``fattened vertices'' of the form $V_\eps:=R_\eps \cdot V$, with some open, bounded set $V\subset\R^N$ and a scale parameter $R_\eps\to 0$ for $\eps\to 0$. This geometric configuration has been studied in \cite{KZ,EP} who proved spectral convergence for the operator $-\Delta$ with Neumann boundary conditions. 
The nature of the limit spectrum depends on the relative scaling of the edge neighbourhoods $E_\eps$ and the vertex neighbourhoods $V_\eps$. 
\begin{enumerate}[(i)]
\item if $|V_\eps|/|E_\eps|\to 0$, the limit spectrum is that of the graph Laplacian with Neumann-Kirchhoff vertex conditions; 
\item if $|V_\eps|/|E_\eps|\to \infty$, the different edges decouple in the limit and the limit spectrum will be the union of the Dirichlet-spectra of all individual edges; 
 \item if $|V_\eps|/|E_\eps|\to q>0$, the spectrum converges to the solution $(u,\lambda)$ of the problem
\begin{equation}\label{zeng}
	\begin{cases}
		u''=\lambda u &\text{\; on\; each\; edge\; }e\\
		\sum_{e\ni v} u'_e(v) = \lambda q u(v), &\text{\; at\; each\; vertex\; }v,
	\end{cases}
\end{equation}
where the sum is over all edges $e$ ending on $v$ and $u'_e(v)=\lim_{x\to v,x\in e}u'(x)$. Since the spectral parameter $\lambda$ appears in the vertex condition, this is a \emph{generalised eigenvalue problem}.
\end{enumerate}
The notion of norm-resolvent convergence in the cases (i), (ii), (iii) has been studied in \cite{P12}

In the following we will apply our above results to study the influence of perforations on fattened graphs.
\subsection{Building the fattened graph}\label{sec:graph-like_Domains1}
Let us first describe in detail how the fattened graph is defined. Let $\Gamma$ be a finite, connected metric graph embedded in $\R^N$. We will give a local description of its fattened analogue around an arbitrary vertex $v\in\Gamma$. Denote by $e_1,\dots, e_{n_v}$ all edges in $\Gamma$ incident to $v$ and let $\ell_1,\dots,\ell_{n_v}$ denote their lengths. Every $e_i$ is canonically isometric to the line segment $\{0\}\times(0,\ell_i)\subset\R^{N-1}\times\R$ via an orthogonal transformation $\Theta_i$ that is unique up to rotation around $e_i$, followed by a shift by $v$.
To build the fattened edges, let $\Om_0$ be as in Section \ref{sec:geometric_setting} with $0\in\Om_0$ and, for every $i\in\{1,\dots,n_v\}$, fix an orthogonal transformation $\Theta_i$ as just described. For $\eps>0$, we call the sets $E_{\eps,i}:=\Theta_i(\eps\Om_0\times(0,\ell_i))+v$ \emph{edge neighbourhoods}. For simplicity we take the same set $\Om_0$ for all edges here.
Similarly, in appropriately shifted coordinates in which $v=0$, we choose a connected, open, bounded set $V\subset \R^N$ with $C^1$ boundary such that $0\in V$.
We call the scaled set $R_\eps V$ a \emph{vertex neighbourhood} of $v$. 
For technical reasons we make the additional assumption that for all $\eps>0$ $\overline V_\eps$ intersects each edge ``only once'', i.e. for all $j\in\{1,\dots,n_v\}$ the implication
\begin{equation}\label{eq:only_once}
	x\in e_j\setminus\overline V_\eps\quad \Rightarrow\quad y\notin \overline V_\eps \;\forall\; y\in e_j \text{ with } |y-v|>|x-v|
\end{equation} 
holds. We note that the set $V$ may be different for every vertex $v\in\Gamma$, while the scaling factor $R_\eps$ is assumed to be global.

In the case $R_\eps\sim\eps$, we make the additional assumption that $\del V$ contains $n_v$ flat copies $\{F_1,\dots,F_{n_v}\}$ of $\Om_0$ such that $F_j\cap e_j = \Theta_j(\Om_0\times\{t\})+v$ for some $t>0$ (these will serve as ``docking sites'' for the edge neighbourhoods). 
In all other cases, where $\nicefrac{\eps}{R_\eps}\to0$, this last assumption on $V$ is unnecessary, since the edge neighbourhoods can be attached to $V_\eps$ via small collars, as the following lemma shows.
\begin{lemma}\label{lemma:bottleneck}
	Let $\nicefrac{\eps}{R_\eps}\to 0$ and let $V_\eps = R_\eps V$, where $V$ is a connected, open, bounded set $V\subset \R^N$ with $C^1$ boundary. 
	If $\eps$ is small enough, then for each edge neighbourhood $E_{\eps,j}$ there exists a $\mathcal O(R_\eps)$ shift $\eta_{\eps,j}\in\R^N$ and a collar domain $B_{\eps,j}$ joining $E_{\eps,j}+\eta_{\eps,j}$ to $V_\eps$ such that $(E_{\eps,j}+\eta_{\eps,j})\cap V_\eps=\emptyset$ and the length $d_{\eps,j}$ of $B_{\eps,j}$ is bounded by
	\begin{align}\label{eq:bottleneck}
			d_{\eps,j} \leq R_\eps\diam(V)
	\end{align}
	(cf. Figure \ref{fig:bottleneck}). In particular $d_{\eps,j} \to 0$ as $\eps\to 0$ for any $j$.
\end{lemma}
\begin{proof}
	Without loss of generality, assume that $R_\eps\diam(V)<\min\{\ell_1,\dots,\ell_{n_v}\}$. Let $\eta_{\eps,j}$ denote the minimiser of the set $\big\{|\eta|\,\big|\,\eta$ parallel to $e_j$ and $\overline V_\eps\cap (E_{\eps,j}+\eta)=\emptyset \big\}$. Then, clearly, $|\eta_{\eps,j}|\leq \diam (V_\eps) = R_\eps\diam(V)$.
	\begin{figure}[h]
		\centering
%
%
%

\begin{tikzpicture}[>=stealth]

	\clip (-1,-3.9) rectangle (\xmax,3.5);
	
	\begin{scope}[rotate=-5]
	\fill[fill=gray!30, rotate=10] (2.15, -1) rectangle (\xmax, 1);
	\fill[fill=gray!80, rotate=10] (1.05, -1) rectangle (2.15, 1);
	
	\begin{scope}[rotate=10]
		\node [] at(\xmax-1,0) {\small$\varepsilon\operatorname{diam}(\Omega_0)$} ;
		\draw[<-] (\xmax-1,1)  -- (\xmax-1,0.3);
		\draw[<-] (\xmax-1,-1) -- (\xmax-1,-0.3);
	\end{scope}
	
	\draw[black, name path=A, line width=1.5pt, domain=-3.5:4, samples=300, shift={(-0.13,0.06)}, rotate=275] 
		plot (\x,{-0.3*sin(pi*30*\x-100) - (\x/4)^5+1.9}) (0, 0);
	\fill[white, name path=B, domain=-4:4, samples=300, shift={(-0.13,0.06)}, rotate=275] 
		plot (\x,-500) (0, 0);
	\tikzfillbetween[of=A and B]{white};
		
	\tikzfillbetween[of=A and B]{gray!20, opacity=0.3};
	
	\node [] at(1.4,3.2) {\small$\partial V_\eps$} ;
	\node [] at(-0.5,3) {\small$V_\eps$} ;
	
	\draw[color=gray, dashed, rotate=10] (-0.5, -1) -- (2.15, -1);
	\draw[color=gray, dashed, rotate=10] (-0.5, 1) -- (2.15, 1);
	
	\draw[black, line width=1.5pt, domain=-3.5:4, samples=300, shift={(-0.13,0.06)}, rotate=275] 
		plot (\x,{-0.3*sin(pi*30*\x-100) - (\x/4)^5+1.9}) (0, 0);
	
	\begin{scope}[rotate=10]
		\node [fill=white] at(2.05,\boxy-0.2) {\small$d_{\eps,j}$} ;
		\draw[->] (2.5,\boxy) -- (2.2,\boxy);
		\draw[dotted] (2.2, \boxy)  -- (2.2,-1);
		
		\draw[->] (1.4,\boxy) -- (1.7,\boxy);
		\draw[dotted] (1.7, \boxy)  -- (1.7,-1);

		\node [fill=white, rotate=5, inner sep=0pt] at(0.9,\boxy-1) {\small$\leq R_\varepsilon\operatorname{diam}(V)$} ;
		\draw[->] (2.6,\boxy-1) -- (2.2,\boxy-1);
		\draw[dotted] (2.2, \boxy-1)  -- (2.2,-1);
		
		\draw[->] (-0.9,\boxy-1) -- (-0.5,\boxy-1);
		\draw[dotted] (-0.5, \boxy-1)  -- (-0.5,-1);
	\end{scope}
	
	\draw[color=black] (1.75, 1.2) -- (3.2,3);
	\node [fill=white] at(3.2,2.8) {\small$B_{\varepsilon,j}$} ;
	\node [] at(\xmax-3,1.5) {\small$E_{\varepsilon, j}$} ;
	
	\node at(-0.2,-0.2) {\small$v$} ;
	\fill[rotate=10] (-0.5,-0.2)circle(.05);

	\end{scope}
	
\end{tikzpicture}

		\caption{Sketch of collar for $\eps\ll R_\eps$.}
		\label{fig:bottleneck}
	\end{figure}
	A collar $B_{\eps,j}$ can now be defined as 
	 $B_{\eps,j}=(\Theta_j(\eps\Om_0\times(0,|\eta_{\eps,j}|))+v)\setminus V_\eps$.
	By construction the length of $B_{\eps,j}$ is bounded by $|\eta_{\eps,j}|$. Finally, note that by our assumptions on $V_\eps$, we have that $(E_{\eps,j}+\eta_{\eps,j}) \cap V_\eps=\emptyset$ for $\eps$ small enough. This follows from \eqref{eq:only_once} and the fact that $\nicefrac{\eps}{R_\eps}\to 0$.
\end{proof}
Similar methods of flattening or attaching collars to the vertex neighbourhoods have been used in the literature (cf. \cite[Sec. 6]{EP}, \cite[Sec. 3.2]{KZ}).
In the following sections, we will assume that such collars $B_{\eps,j}$ are used to define the fattened graph whenever $\nicefrac{\eps}{R_\eps} \to 0$. To streamline notation, we define $B_{\eps,j}:=\emptyset\;\forall j$ when $R_\eps\sim\eps$.
\begin{de}
Given a finite, connected graph $\Gamma$, by a \emph{fattened analogue} we shall mean a family of open subsets of $\R^N$ (indexed by $\eps>0$), consisting of edge neighbourhoods $E_{\eps,j}$ and vertex neighbourhoods $V_\eps$, which are linked according to the connection rules of $\Gamma$, using the techniques described above. For every edge $E_j$, there will be two collars, $B_{\eps,j}^l$ (attached at $\eps\Om_0\times\{0\}$) and $B_{\eps,j}^r$ (attached at $\eps\Om_0\times\{\ell_j\}$).
\end{de}
\begin{remark}
	\begin{enumerate}[(i)]
		\item
		According to Lemma \ref{lemma:bottleneck}, the fattened edges and vertices have to be slightly moved with respect to their original counterparts. We will ignore this in our notation in the following, since all equations considered are invariant under shifts. That is, instead of $E_{\eps,j}+\eta_{\eps,j}$ we simply write $E_{\eps,j}$, etc.
		\item
		When building the fattened graph via Lemma \ref{lemma:bottleneck}, the shifts $\eta_{\eps,j}$ will in general change the angles between the edges. This does not affect the results in the following sections, because the graph Laplacians defined in eqs. \eqref{eq:Fast_graph_limit}, \eqref{eq:slow_graph_limit} and \eqref{eq:borderline_graph_limit} depend only on the metric graph structure of $\Gamma$ (that is, the connection rules and the lengths of the edges) and are independent of the particular embedding in $\R^N$.
	\end{enumerate}
\end{remark}
\subsection{Small vertex neighbourhoods}\label{sec:small}
Let us first consider the situation in which $\nicefrac{|V_\eps|}{|E_\eps|}\to 0$. To be precise, we assume in this section that 
\begin{align*}
	\eps \leq R_\eps = o\big(\eps^{\f{N-1}{N}}\big).
\end{align*}
The lower bound on $R_\eps$ ensures that the diameter of $V_\eps$ scales at least as the diameter of the $E_{\eps,j}$, i.e. the edge neighbourhoods do not overlap as $\eps\to 0$.

Let $\Gamma$ be a finite, connected metric graph and denote by $\Om_\eps$ a fattened analogue. Let $v$ be a vertex of $\Gamma$ and $e_1,\dots, e_n$ be all edges incident to $v$ with lengths $\ell_1,\dots, \ell_n$.
\begin{figure}[htbp]
	\centering
%

\begin{tikzpicture}[>=stealth]
	\clip (-10.5,-4) rectangle (5.9,4);
	\begin{scope}[xshift=10]
		\draw (-9,0) -- (-6.3,0);
		\node [] at(-6,0) {$\cdots$} ;
		\draw (-9,0) -- (-11.5,3.5);
		\draw (-9,0) -- (-11.5,-3.5);
		
		\fill (-9,0) circle(1mm);
		
		\node [] at(-8.8,0.2) {$v$} ;
		\node [] at(-7,0.2) {$e_1$} ;
		\node [] at(-10.5,2.6) {$e_2$} ;
		\node [] at(-10.5,-2.6) {$e_3$} ;
	\end{scope}
	
	\fill[fill=gray!30, rotate=10] (1, -1) rectangle (7, 1);
	\fill[fill=gray!30,rotate=240] (1,-1) rectangle (6,1);
	\fill[fill=gray!30,rotate=120] (1,-1) rectangle (6,1);
	
	\fill[fill=gray!90, rotate=10] (1, -1) rectangle (2.6, 1);
	\fill[fill=gray!90, rotate=240] (1, -1) rectangle (2.7, 1);
	\fill[fill=gray!90, rotate=120] (1, -1) rectangle (2.5, 1);
	
	\fill[fill=gray!50] plot [smooth cycle, tension=0.4, xshift=-8, yshift=-3] coordinates {(-2.5,0) (-1,2) (1.2,2.3) (2.2,2) (2.5,0) (1.5,-2) (0,-2) (-1,-2) };
	
	\foreach \y in { -0.5, 0.5, 0} {
	\foreach \x in {2, 2.5, 3, 3.5, 4, 4.5} {
		\draw[xshift=30, yshift=4, rotate=10] (\x,\y) circle (.5mm);
		\fill[color=white, xshift=30, yshift=4,  rotate=10] (\x,\y) circle (.5mm);
	}}
	\foreach \y in { -0.5, 0.5, 0} {
	\foreach \x in {2.5, 3, 3.5, 4, 4.5} {
		\draw[rotate=240, xshift=20] (\x,\y) circle (.5mm);
		\fill[color=white, rotate=240, xshift=20] (\x,\y) circle (.5mm);
	}}
	\foreach \y in {-0.5, 0.5, 0} {
	\foreach \x in {2.5, 3, 3.5, 4, 4.5, 5, 5.5, 6} {
		\draw[rotate=120,xshift=10] (\x,\y) circle (.5mm);
		\fill[color=white, rotate=120, xshift=10] (\x,\y) circle (.5mm);
	}}
	
	\foreach \x in {-2,-1.5,-1,-0.5,0,0.5,1,1.5} {
		\foreach \y in {-0.5, 0, 0.5}{
			\draw (\x,\y) circle (.5mm);
			\fill[color=white] (\x,\y) circle (.5mm);
		}
	}
	\foreach \x in {-1.5, -1, -0.5, 0, 0.5, 1, 1.5} {
		\foreach \y in {-1, 1}{
			\draw (\x,\y) circle (.5mm);
			\fill[color=white] (\x,\y) circle (.5mm);
		}
	}
	\foreach \x in {-1, -0.5, 0, 0.5, 1} {
		\foreach \y in {-1.5, 1.5}{
			\draw (\x,\y) circle (.5mm);
			\fill[color=white] (\x,\y) circle (.5mm);
		}
	}

	\node [fill=white,inner sep=1pt] at(-3.15,3.45) {$E_{\varepsilon,1}$} ;
	\node [fill=white,inner sep=1pt] at(4.1,1) {$E_{\varepsilon,2}$} ;
	\node [fill=white,inner sep=1pt] at(-3.15,-3.45) {$E_{\varepsilon,3}$} ;
	\node [fill=white,inner sep=1pt] at(0,0) {$V_\varepsilon$} ;
	
	\node [fill=white,inner sep=1pt] at(3,3) {$B_{\varepsilon,2}^l$} ;
	\node [fill=white,inner sep=1pt] at(1,-3) {$B_{\varepsilon,3}^r$} ;
	\draw (2.3,0.9) -- (3,2.7);
	\draw (-0.5,-2.5) -- (0.6,-2.9);
	
	\draw[<->] (5.1,1.9) -- (5.45,-0.05);
	\node [fill=white] at(5.3,0.9) {$\sim \varepsilon$} ;
	
	\draw[dotted, xshift=10] (-1,-2.2) -- (-4,-2.2);
	\draw[dotted, xshift=10] (0,2.2) -- (-4,2.2);
	\draw[<->,    xshift=10] (-4,-2.2) -- (-4, 2.2);
	\node[fill=white, xshift=10] at(-4,0) {$\sim R_\varepsilon$} ;

\end{tikzpicture}

	\caption{Sketch of graph like perforated domain. The relative scaling between $R_\eps$ and $\eps$ is different in each subsection.}
	\label{fig:small}
\end{figure}
As discussed in Section \ref{sec:graph-like_Domains1}, after suitable changes of coordinates the vertex neighbourhood is of the form $V_\eps = R_\eps\!\cdot\! V$ with $\f{R_\eps^N}{\eps^{N-1}}\to 0$ as $\eps\to 0$ and the fattened edges are of the form $E_{\eps,i} = (\eps\Om_0)\times (0,\ell_i)$. Introducing a periodic perforation $T_\eps$ as shown in Figure \ref{fig:small} defines a domain $\Om_\eps^{\text p}$.
\begin{remark}
	On each edge neighbourhood we choose the perforation to be aligned with the corresponding edge, in order to be able to apply the results of Section \ref{sec:proof_of_theorem}. The perforation of the vertex neighbourhood can be chosen with arbitrary orientation without affecting the limit. This follows from the fact that the classical homogenisation results hold for arbitrary domains (cf. \cite{CM}).
\end{remark}
 Note that we do not perforate the collars $B_{\eps,j}^{l,r}$. On this domain we consider the Poisson equation with Dirichlet boundary conditions on the holes.
\begin{align}\label{eq:graph_problem}
	\begin{cases}
		(-\Delta+z) u_\eps = f_\eps & \text{ in }\Oeps\\
		\hfill u_\eps = 0 & \text{ on } \del T_\eps\\
		\hfill \del_\nu u_\eps = 0 &\text{ on } \del\Om_\eps
	\end{cases}
\end{align}
for $z>0$ and $f_\eps\in L^2(\Om_\eps)$ with $\|f_\eps\|_{L^2(\Om_\eps)}$ uniformly bounded. 

This new geometric situation requires new identification operators to be defined. To this end, let $L^2(\Gamma):=\bigoplus_{j=1}^{n_{\mathrm{e}}} L^2(e_j)$, where $\{e_j\}_{j=1}^{n_{\mathrm{e}}}$ is the set of edges of $\Gamma$, and let $H^1(\Gamma)$ denote the space of continuous functions $\phi$ on $\Gamma$ such that for every edge $e_j$ the restriction $\phi|_{e_j}$ is in $H^1(e_j)$. Moreover, let us define
\begin{equation}\label{eq:UGammaDef}
\begin{split}
	&\U^\Gamma : L^2(\Gamma) \to L^2(\Om_\eps)\\
	&\U^\Gamma\phi(x) = |\eps\Om_0|^{-\f12}\cdot\begin{cases}
		\phi(t) & \text{ if } x=(\bar x,t)\in E_{\eps,j},\; {j\in\{1,\dots,n_{\mathrm{e}}\}}\\
		0 & \text{ if } x\in V_\eps {\cup \bigcup_{\{j:e_j\ni v\}}B_{\eps,j}^{l,r},}
	\end{cases}
\end{split}
\end{equation}
where $(\overline x,t)$ are understood to mean local coordinates running along the fattened edge, that is $\overline x\in \eps\Om_0$, $t\in(0,\ell_j)$, as described in Section \ref{sec:graph-like_Domains1}.
In the union $\bigcup_{e_j\ni v}B_{\eps,j}^{l,r}$ we include either $B_{\eps,j}^{l}$ or $B_{\eps,j}^{r}$, depending on which end of $e_j$ meets $v$. In other words, the union is over all collars that meet $V_\eps$.
Problem \eqref{eq:graph_problem} immediately yields the a priori bound
\begin{align}\label{eq:apriori_fast}
	\|\nabla u_\eps\|^2_{L^2(\Om_\eps)} \leq C \|f_\eps\|^2_{L^2(\Om_\eps)}.
\end{align}
A proof analogous to that of Proposition \ref{weakconv} shows that there exists a subsequence (again denoted by $u_\eps$) such that $\| u_\eps-\U^\Gamma u\|_{L^2(\Om_\eps)}\to 0$ for some $u\in H^1(\Gamma)$. Note that the fact that $|V_\eps|/|E_\eps|\to 0$ ensures the convergence on the vertex neighbourhoods.

We are now going to derive an equation on $\Gamma$ that identifies the limit $u$. To this end, we define a second identification operator $\mathcal V_\eps^\Gamma$ which preserves $H^1$ regularity. Let
\begin{align*}
	&\V^\Gamma : H^1(\Gamma) \to H^1(\Om_\eps)\\
	&\V^\Gamma\phi(x) = |\eps\Om_0|^{-\f12}\cdot\begin{cases}
		\phi(t) & \text{ if } x=(\bar x,t)\in E_{\eps,j},\;j\in\{1,\dots,n_{\mathrm{e}}\}\\
		\phi(v) & \text{ if } x\in V_\eps {\cup \bigcup_{\{j:e_j\ni v\}}B_{\eps,j}^{l,r}},
	\end{cases}
\end{align*}
Let $w_\eps$ now be defined as in \eqref{wepsilon} ($w_\eps\equiv 1$ on the $B_{\eps,j}^{l,r}$) and consider the weak formulation of this problem with test function $w_\eps \V^\Gamma\phi$ for arbitrary $\phi\in H^1(\Gamma)$. Note that $w_\eps \V^\Gamma\phi\in H^1(\Om_\eps)$ with $w_\eps \V^\Gamma\phi=0$ on the holes, and is therefore a valid test function for the perforated domain problem. The weak formulation of \eqref{eq:graph_problem} now reads
\begin{align*}
	\int_{\Om_\eps^{\text p}} \overline{\nabla u_\eps}\cdot \nabla\big(w_\eps\V^\Gamma\phi\big)\,dx &= \int_{\Om_\eps^{\text p}} \overline f_\eps w_\eps\V^\Gamma\phi \,dx
\end{align*}
Decomposing into the different components of $\Om_\eps^{\text p}$ we obtain
\begin{equation}\label{eq:small_weak_form}
\begin{aligned}
	&\sum_{i=1}^{n_{\mathrm{e}}}\int_{E_{\eps,i}} \overline{\nabla u_\eps} \cdot\nabla \big(w_\eps\V^\Gamma\phi\big)\,dx + 
	\sum_{i=1}^{n_{\mathrm{e}}}\int_{B_{\eps,i}^l\cup B_{\eps,i}^r} \overline{\nabla u_\eps} \cdot\nabla \big(w_\eps\V^\Gamma\phi\big)\,dx + 
	\sum_{j=1}^{n_\mathrm{v}}\int_{V_{\eps,j}} \overline{\nabla u_\eps}\cdot\nabla \big(w_\eps\V^\Gamma\phi\big)\,dx\\ 
	&+ z\sum_{i=1}^{n_{\mathrm{e}}}\int_{E_{\eps,i}} \overline u_\eps w_\eps\V^\Gamma\phi \,dx + 
	z\sum_{i=1}^{n_{\mathrm{e}}}\int_{B_{\eps,i}^l\cup B_{\eps,i}^r} \overline u_\eps w_\eps\V^\Gamma\phi \,dx  + 
	z\sum_{j=1}^{n_\mathrm{v}}\int_{V_{\eps,j}} \overline u_\eps w_\eps\V^\Gamma\phi\,dx\\
	&= \sum_{i=1}^{n_{\mathrm{e}}}\int_{E_{\eps,i}} \overline f_\eps w_\eps\V^\Gamma\phi \,dx + 
	\sum_{i=1}^{n_{\mathrm{e}}}\int_{B_{\eps,i}^l\cup B_{\eps,i}^r} \overline f_\eps w_\eps\V^\Gamma\phi \,dx + 
	\sum_{j=1}^{n_\mathrm{v}}\int_{V_{\eps,j}} \overline f_\eps w_\eps\V^\Gamma\phi\,dx
\end{aligned}
\end{equation}
for all $\phi\in H^1(\Gamma)$, where $n_{\mathrm{e}},\,n_\mathrm{v}$ denote the number of edges and vertices of $\Gamma$, respectively. 
Let us next show that all integrals over the collars $B_{\eps,i}^l\cup B_{\eps,i}^r$ do not contribute to the limit.
First, note that all the terms $\int_{B_{\eps,i}^{l,r}} \overline{\nabla u_\eps} \cdot\nabla \big(w_\eps\V^\Gamma\phi\big)\,dx$ vanish identically, because $w_\eps\V^\Gamma\phi$ is constant on $B_{\eps,i}^{l,r}$. Moreover, the terms
\begin{align*}
	z\sum_{i=1}^{n_{\mathrm{e}}}\int_{B_{\eps,i}^{l,r}} \overline u_\eps w_\eps\V^\Gamma\phi \,dx 
\end{align*}
from the second line of \eqref{eq:small_weak_form} can be estimated as follows.
\begin{align*}
	\left|\int_{B_{\eps,i}^{l,r}} \overline u_\eps w_\eps\V^\Gamma\phi \,dx \right|
	&\leq \|u_\eps\|_{L^2(B_{\eps,i}^{l,r})}\|\V^\Gamma\phi\|_{L^2(B_{\eps,i}^{l,r})} \\
	&= \|u_\eps\|_{L^2(B_{\eps,i}^{l,r})}|\eps\Om_0|^{-\f12}\|\phi(v)\|_{L^2(B_{\eps,i}^{l,r})} \\
	&\leq \|u_\eps\|_{L^2(B_{\eps,i}^{l,r})}|\eps\Om_0|^{-\f12} |\phi(v)|\, \big|B_{\eps,i}^{l,r}\big|^{\f12},
\end{align*}
where we have used the fact that $w_\eps\equiv 1$ on $B_{\eps,i}^{l,r}$ in the first line. Note that the measure $\big|B_{\eps,i}^{l,r}\big|$ is equal to $|\eps\Om_0|\cdot d_{\eps,j}$ (recall the definition of $d_{\eps,j}$ from Lemma \ref{lemma:bottleneck}). Thus, we get
\begin{align*}
	\left|\int_{B_{\eps,i}^{l,r}} \overline u_\eps w_\eps\V^\Gamma\phi \,dx \right|
	&\leq \|u_\eps\|_{L^2(B_{\eps,i}^{l,r})}|\eps\Om_0|^{-\f12} |\phi(v)| |\eps\Om_0|^{\f12}\cdot d_{\eps,j}^{\f12}\\
	&= \|u_\eps\|_{L^2(B_{\eps,i}^{l,r})}|\phi(v)| d_{\eps,j}^{\f12}.
\end{align*}
Since $d_{\eps,j}\to 0$ as $\eps\to 0$, by Lemma \ref{lemma:bottleneck}, and $\|u_\eps\|_{L^2(B_{\eps,i}^{l,r})}$ is bounded, we conclude that $\int_{B_{\eps,i}^{l,r}} \overline u_\eps w_\eps\V^\Gamma\phi \,dx\to 0$ for all $i$ as $\eps\to 0$. 
An analogous argument shows that $\sum_i\int_{B_{\eps,i}^l\cup B_{\eps,i}^r} \overline f_\eps w_\eps\V^\Gamma\phi \,dx \to 0$ as $\eps\to 0$.

Next we turn to the integrals over the $E_{\eps,i}$ and $V_{\eps,j}$.
Since every fattened edge is of the form $E_{\eps,i} = (\eps\Om_0)\times (0,\ell_i)$, we can immediately conclude from the proof of Theorem \ref{normconv} that 
\begin{equation}\label{eq:fast_edge_convergences}
\begin{split}
	\sum_{i=1}^{n_{\mathrm{e}}}\int_{E_{i,\eps}} \overline{\nabla u_\eps}\cdot\nabla \big(w_\eps\V^\Gamma\phi\big)\,dx &\,\to\, \sum_{i=1}^{n_{\mathrm{e}}}\int_{e_i} \overline{\nabla u}\cdot\nabla \phi\,dt + \mu\sum_{i=1}^{n_{\mathrm{e}}}\int_{e_i} \overline{u} \phi\,dt \qquad\text{and}\\
	\sum_{i=1}^{n_{\mathrm{e}}}\int_{E_{i,\eps}} \overline f_\eps w_\eps\V^\Gamma\phi \,dx &\,\to\, \sum_{i=1}^{n_{\mathrm{e}}}\int_{e_i} \overline f\phi\,dt\\
	z\sum_{i=1}^{n_{\mathrm{e}}}\int_{E_{i,\eps}} \overline u_\eps w_\eps\V^\Gamma\phi \,dx &\,\to\, z\sum_{i=1}^{n_{\mathrm{e}}}\int_{e_i} \overline u\phi\,dt
\end{split}
\end{equation}
whenever $f_\eps\xrightharpoonup{L^2}f$ on each edge.  It remains to study the integrals over $V_{\eps,j}$. To treat the gradient term, let $j\in\{1,\dots,n_\mathrm{v}\}$ and compute
\begin{align*}
	\left| \int_{V_{\eps,j}} \overline{\nabla u_\eps}\cdot\nabla \big(w_\eps\V^\Gamma\phi\big)\,dx \right| &= \left|\int_{{V_\eps,j}} \overline{\nabla u_\eps}\cdot\nabla w_\eps\big(\V^\Gamma\phi\big) \,dx + \int_{V_{\eps,j}} \overline{\nabla u_\eps}\cdot\nabla\big(\V^\Gamma\phi\big)\, w_\eps\,dx\right|\\
	&= \left|\int_{V_{\eps,j}} \overline{\nabla u_\eps}\cdot\nabla w_\eps\big(\V^\Gamma\phi\big) \,dx\right|\\
	&\leq C\|\nabla u_\eps \|_{L^2(V_{\eps,j})} \bigl\|\eps^{\f{-N+1}{2}}\nabla w_\eps\bigr\|_{L^2(V_{\eps,j})} |\phi(v)| \\
	&\leq C \|f_\eps\|^2_{L^2(\Om_\eps)}\bigl\|\eps^{\f{-N+1}{2}}\nabla w_\eps\bigr\|_{L^2(V_{\eps,j})} |\phi(v)|\\
	&\leq C \bigl\|\eps^{\f{-N+1}{2}}\nabla w_\eps\bigr\|_{L^2(V_{\eps,j})},
\end{align*}
where we have used \eqref{eq:apriori_fast} in the fourth line. An explicit computation shows that 
\begin{align*}
	\big\|\eps^{\f{-N+1}{2}}\nabla w_\eps\big\|_{L^2(V_{\eps,j})}^2 \leq C\f{R_\eps^N}{\eps^{N-1}}.
\end{align*}
 Thus, the term $\int_{V_{\eps,j}} \overline{\nabla u_\eps}\nabla \big(w_\eps\V^\Gamma\phi\big)\,dx$ converges to 0 as $\eps\to 0$. Similarly, we compute
\begin{align*}
	\int_{V_{\eps,j}} \overline f_\eps w_\eps\, \V^\Gamma\phi\,dx &\leq \|f_\eps\|_{L^2(\Om_\eps)}|\phi(v)|\eps^{\f{-N+1}{2}}\|w_\eps\|_{L^2(V_{\eps,j})}\\
	&\leq C\eps^{\f{-N+1}{2}}|V_{\eps,j}|^{\f12}\\
	&\to 0
\end{align*}
as $\eps\to 0$. Finally, we have
\begin{align*}
	z\left|\int_{V_{\eps,j}} \overline u_\eps w_\eps\V^\Gamma\phi\,dx\right| &\leq z\|f_\eps\|_{L^2(\Om_\eps)}|\phi(v)|\eps^{\f{-N+1}{2}}\|w_\eps\|_{L^2(V_{\eps,j})}\\
	&\leq zC\eps^{\f{-N+1}{2}}|V_{\eps,j}|^{\f12}\\
	&\to 0
\end{align*}
as $\eps\to 0$. 
Since the vertex $v_j$ was arbitrary in the above procedure, we conclude that the limit $u$ solves the problem
\begin{align}\label{eq:Fast_graph_limit}
	\int_{\Gamma} \overline{\nabla u}\nabla \phi\,dt + (z+\mu) \int_{\Gamma} \overline{u} \phi\,dt =  \int_{\Gamma} \overline f\phi\,dt\qquad\forall\phi\in H^1(\Gamma),
\end{align}
which is nothing but the sesquilinear form of the operator $-\Delta+\mu$ on $L^2(\Gamma)$ with Neumann-Kirchhoff boundary conditions at each vertex. Since we only used weak $L^2$-convergence of $f_\eps$, we can argue as in the proof of Lemma \ref{uniformU} to obtain a norm-resolvent convergence statement. More precisely, if we define
\begin{equation}\label{eq:Adef_fast}
\begin{aligned}
	A_\eps^\Gamma & :=-\Delta,  &\;\mathcal D(A_\eps^\Gamma)&=\big\{u\in H^2(\Oeps) : \del_\nu u|_{\del\Om_\eps}=0 \text{ and } u|_{\del T_\eps}=0\big\}\\
	A^\Gamma  & := -\Delta+\mu, &\;\mathcal D(A^\Gamma)&=\Big\{ u\in H^2(\Gamma) :  \sum_{e\ni v} u'_e(v)=0 \text{ at all vertices }v\Big\},
\end{aligned}
\end{equation}
(where $H^2(\Gamma)$ is a defined as $C(\Gamma)\cap \bigoplus_{i=1}^{n_{\mathrm{e}}} H^2(e_i)$)
then we have the following
\begin{theorem}
	If $\f{R_\eps^N}{\eps^{N-1}}\to 0$ as $\eps \to 0$, then
	\begin{align*}
		\left\| (A_\eps^\Gamma+z)^{-1}\U^\Gamma - \U^\Gamma (A^\Gamma+z)^{-1} \right\|_{\mathcal L(L^2(\Gamma),L^2(\Oeps))} &\to 0 
	\end{align*}
	as $\eps\to 0$.
\end{theorem}
It is easily seen that the conditions for Theorem \ref{IOSThm} are also satisfied by the pair $(A_\eps^\Gamma,\U^\Gamma)$, which allows us to conclude that
\begin{corollary}\label{SpCon_Gamma}
	Choose $z=1$ and let $\lambda_k^\eps$ and $\lambda_k$ denote the $k$-th eigenvalues of $A_\eps^\Gamma$ and $A^\Gamma$, respectively. There exist a constant $C>0$ and a function $a(\eps)$ with $a(\eps)\to 0$ as $\eps\to 0$ such that 
	\begin{align*}
		|(\lambda_k^\eps)^{-1} - \lambda_k^{-1}| \leq C a(\eps)\qquad\text{ for all }k\in\mathbb N,
	\end{align*}
	where $C$ is independent of $\eps$ and $k$.
\end{corollary}
\subsection{Large vertex neighbourhoods}\label{sec:large_vertex}
Next, we study the case of large vertex neighbourhoods, i.e. $|V_\eps|/|E_\eps|\to\infty$. In other words, we assume $V_\eps = R_\eps\cdot V$ for some open, bounded set $V$ as in Section \ref{sec:graph-like_Domains1}, where  $\f{R_\eps^N}{\eps^{N-1}}\to \infty$ as $\eps\to 0$. Here the situation is different from that in the previous subsection because the vertex neighbourhoods cannot be neglected in the limit anymore. In particular, spectral convergence will not {follow straightforwardly in this case}, since $(\U^\Gamma)$ does not satisfy (H4) in Theorem \ref{IOSThm} for large vertex neighbourhoods. 
Indeed, spectral convergence in a narrow sense is expected to fail, as this is already the case in the classical situation (without perforation).
 This is easily seen from the fact that the Neumann Laplacians on the graph like domain all have 0 as an eigenvalue, whereas the limit operator (a decoupled \emph{Dirichlet} Laplacian) does not. 
In the classical case this fact is circumvented by considering dilated versions of the operators involved in order to re-introduce the 0 eigenvalue on the graph (see for instance \cite[Sec. 6, 7]{EP}). The question to what extent those methods can be applied to the perforated case, will be studied in future work, but here we shall content ourselves with proving only strong convergence. Similar comments apply to the borderline case which is studied in the next section.
To prove strong convergence, let $f\in L^2(\Gamma)$ and consider the equation 
\begin{align}\label{eq:large_equation}
	(A_\eps+z)u_\eps=\U^\Gamma f
\end{align}
on $\Om_\eps$. As a preparation, note that from the a priori estimate \eqref{eq:apriori_fast} we obtain a bound for $u_\eps$ on the vertex neighbourhoods 
\begin{align}
	\|\nabla u_\eps\|_{L^2(V_\eps)} \leq C \|f\|_{L^2(\Gamma)}.
\end{align}
A blow up argument as in the proof of Proposition \ref{weakconv} shows that {for any vertex $v$} there exists a constant $u_v$ such that $\big\|u_\eps-|V_\eps|^{-\nicefrac{1}{2}}u_v\big\|_{L^2(V_\eps)}\to 0$. We will show that necessarily $u_v=0$. Owing to the new scale $|V_\eps|$ present in this case, we introduce the extension operator
\begin{equation}\label{eq:W_Gamma_def}
\begin{split}
	&\W^\Gamma : H^1(\Gamma) \to H^1(\Om_\eps)\\
	&\W^\Gamma\phi(x) = |V_\eps|^{-\f12}\cdot\begin{cases}
		\phi(t) & \text{ if } x=(\bar x,t)\in E_{\eps,j},\; j\in\{1,\dots,n_{\mathrm{e}}\}\\
		\phi(v) & \text{ if } x\in V_\eps {\cup \bigcup_{\{j:e_j\ni v\}}B_{\eps,j}^{l,r},}
	\end{cases}
\end{split}
\end{equation}
	where the same comments as below eq. \eqref{eq:UGammaDef} apply to the union $\bigcup_{\{j:e_j\ni v\}}B_{\eps,j}^{l,r}$ and the notation $(\bar x,t)\in E_{\eps,j}$.
To this end, let $\phi\in H^1(\Gamma)$ and $z\neq -\mu$ and use $w_\eps\W^\Gamma\phi$ as a test function in the weak formulation of \eqref{eq:large_equation}.
\begin{align}
	\int_{\Om_\eps} \nabla u_\eps\cdot\nabla(w_\eps\W^\Gamma\phi)\,dx + z\int_{\Om_\eps}u_\eps w_\eps\W^\Gamma\phi\,dx &= \int_{\Om_\eps} (\U^\Gamma f) w_\eps (\W^\Gamma\phi)\, dx \nonumber \\
	&= \sum_{i=1}^{n_{\mathrm{e}}} \int_{E_{i,\eps}} (\U^\Gamma f) w_\eps (\W^\Gamma\phi)\, dx,\label{eq:weak_formulation_large}
\end{align}
where in the last line we used the fact that $\U^\Gamma f=0$ on $V_\eps {\cup \bigcup_{\{j:e_j\ni v\}}B_{\eps,j}^{l,r}} $. As in Lemmas \ref{w=1} and \ref{lemma:laplace_w} one shows that for any $j\in\{1,\dots,n_\mathrm{v}\}$
\begin{align*}
	\int_{V_{\eps,j}} \nabla u_\eps\cdot\nabla(w_\eps\W^\Gamma\phi)\,dx &\to \mu u_{v_j}\phi(v_j)\\
	z\int_{\Om_\eps}u_\eps w_\eps\W^\Gamma\phi\,dx &\to z u_{v_j}\phi(v_j).
\end{align*}
Moreover, all integrals over the edge neighbourhoods $E_{i,\eps}$ converge to 0 by our choice of scaling in \eqref{eq:W_Gamma_def}. 
Similarly, the integrals over the collars $B_{\eps,i}^{l,r}$ vanish in the limit by a similar calculation to that after \eqref{eq:small_weak_form} (with $|\eps\Om_0|^{-\f12}$ replaced by $|V_{\eps,j}|^{-\f12}$), using again Lemma \ref{lemma:bottleneck}.
Therefore, passing to the limit in \eqref{eq:weak_formulation_large} leads to 
\begin{align}\label{eq:mu=z}
	\mu u_v\phi(v) + z u_v\phi(v) = 0 \quad{\text{for any vertex }v\in\Gamma.}
\end{align}
Since $\phi\in H^1(\Gamma)$ was chosen arbitrary and $z\neq\mu$ we conclude from \eqref{eq:mu=z} that $u_v=0$ for all vertices $v$.

Moving on to identifying the limiting equation, we note that it follows from the a priori estimate \eqref{eq:apriori_fast} that on each edge (a subsequence of) $u_\eps\!\!\upharpoonleft_{E_{i,\eps}}$ converges to a function in $H^1(e_i)$. We conclude that there exists a function $u\in \bigoplus_i H^1(e_i)$ such that $\|u_\eps- \U^\Gamma u\|_{L^2(\Om_\eps)}\to 0$. To conclude, we note that since $\|\nabla u_\eps\|_{L^2(\Om_\eps)}$ is uniformly bounded and $u_\eps\to 0$ at each vertex, we must have $u\!\upharpoonleft_{E_{i,\eps}}\in H^1_0(E_{i,\eps})$ for all $i$. 

Finally, we identify the limit equation by letting $\phi\in H^1_0(\Gamma)$ and using $w_\eps\V^\Gamma\phi$ as a test function in the weak formulation of \eqref{eq:large_equation} to obtain
\begin{align}\label{eq:large_weak_formulation}
	\int_{\Om_\eps} \nabla u_\eps\cdot\nabla(w_\eps\V^\Gamma\phi)\,dx + z\int_{\Om_\eps}u_\eps w_\eps\V^\Gamma\phi\,dx &= \int_{\Om_\eps} (\U^\Gamma f) w_\eps (\V^\Gamma\phi)\, dx
\end{align}
By the choice of $\phi$, all integrals over vertex neighbourhoods and collars are zero, while the integrals over the edge neighbourhoods are treated exactly as in the case of small vertex neighbourhoods (cf. \eqref{eq:fast_edge_convergences}). Passing to the limit in \eqref{eq:large_weak_formulation} we conclude that
\begin{align*}
	\int_{\Gamma} \overline{\nabla u}\nabla \phi\,dt + (z+\mu) \int_{\Gamma} \overline{u} \phi\,dt =  \int_{\Gamma} \overline f\phi\,dt\qquad\forall\phi\in \bigoplus_{e\in\Gamma}H^1_0(e).
\end{align*}
To summarise, we have shown that 
\begin{theorem}
	If $\f{R_\eps^N}{\eps^{N-1}}\to\infty$, then for every $f\in L^2(\Gamma)$ one has
	\begin{align*}
		\left\| u_\eps - \U^\Gamma u \right\|_{L^2(\Om_\eps)}\to 0
	\end{align*}
	as $\eps\to 0$, where $u_\eps$ denotes the solution of \eqref{eq:large_equation} and $u\in\bigoplus_{e\in\Gamma}H^1_0(e)$ denotes the solution to the decoupled family of Dirichlet problems 
	\begin{align}\label{eq:slow_graph_limit}
	\begin{cases}
		(-\Delta+\mu+z)u = f &\text{ on e}\\
		\hfill u=0 &\text{ on }\del e
	\end{cases}
	\end{align}
	for all edges $e\in\Gamma$.
\end{theorem}
\subsection{The borderline case ${\nicefrac{|V_\eps|}{|E_\eps|}\to c>0}$}\label{sec:borderline}
Let us now study the case in which the volume of the edge- and the vertex neighbourhoods decay at the same rate. In other words, we assume $V_\eps = R_\eps\cdot V$ for some open, bounded set $V$ as in Section \ref{sec:graph-like_Domains1}, where without loss of generality $\f{R_\eps^N}{\eps^{N-1}}\to 1$ as $\eps\to 0$. We study again problem \eqref{eq:graph_problem} on the corresponding perforated domain.

The discussion before eq. \eqref{eq:apriori_fast} carries over verbatim to the present situation and it only remains to study the integrals over the vertex neighbourhoods and collars. As in Section \ref{sec:small}, we have 
\begin{align}
	\int_{V_{\eps}} \overline{\nabla u_\eps}\cdot\nabla \big(w_\eps\V^\Gamma\phi\big)\,dx  &= \int_{V_\eps} \overline{\nabla u_\eps}\cdot\nabla w_\eps\big(\V^\Gamma\phi\big) \,dx + \int_{V_\eps} \overline{\nabla u_\eps}\cdot\nabla\big(\V^\Gamma\phi\big)\, w_\eps\,dx \nonumber \\
	&= \int_{V_\eps} \overline{\nabla u_\eps}\cdot\nabla w_\eps\big(\V^\Gamma\phi\big) \,dx \label{eq:borderline_weak_formulation}
\end{align}
for any fattened vertex $V_\eps$ and 
\begin{align}
	\sum_{i=1}^{n_{\mathrm{e}}}\int_{B_{\eps,i}^l\cup B_{\eps,i}^r} \overline{\nabla u_\eps}\cdot\nabla \big(w_\eps\V^\Gamma\phi\big)\,dx &= 0
\end{align}
(since $\V^\Gamma\phi$ is constant on $V_\eps$ and $w_\eps\equiv 1$ on the $B_{\eps,i}^{l,r}$), whereas now the right-hand side of \eqref{eq:borderline_weak_formulation} does not converge to zero. As noted in the discussion around eq. \eqref{zeng}, the spectral parameter enters the boundary condition in this case. Hence, the limit operator is not the resolvent of an operator on $L^2(\Gamma)$ and the notion of norm-resolvent convergence makes no sense. Therefore, as in the last subsection, we shall content ourselves with proving strong convergence here. This is readily obtained as follows. The proof of Lemma \ref{lemma:laplace_w} immediately implies that 
\begin{align*}
	\int_{V_\eps} \overline{\nabla u_\eps}\cdot\nabla w_\eps\big(\V^\Gamma\phi\big) \,dx \to \f{|V|}{|\Om_0|}\mu\, \overline u(v)\phi(v)
\end{align*}
for any vertex neighbourhood $V_\eps$.
Finally, we have 
\begin{align*}
	z\int_{V_{\eps}} \overline u_\eps w_\eps\V^\Gamma\phi\,dx \,dx \to \f{|V|}{|\Om_0|}z\, \overline u(v)\phi(v).
\end{align*}
This follows from the facts that $\big\|u_\eps - \V^\Gamma u\big\|_{L^2(V_\eps)}\to 0$ and $\big\|w_\eps\V^\Gamma\phi - \V^\Gamma \phi\big\|_{L^2(V_\eps)}\to 0$. Since $|V_\eps|\sim|E_{i,\eps}|$,  the proofs are entirely analogous to those in Section \ref{sec:sonvergence_of_solutions}.
Hence the weak limit $u$ satisfies the equation
\begin{align}\label{eq:borderline_graph_limit}
	\int_{\Gamma} \overline{\nabla u}\nabla \phi\,dt + (z+\mu) \int_{\Gamma} \overline{u} \phi\,dt + (z+\mu) \f{|V|}{|\Om_0|}\overline u(v)\phi(v)=  \int_{\Gamma} \overline f\phi\,dt\qquad\forall\phi\in H^1(\Gamma),
\end{align}
This is nothing but the sesquilinear form for the Laplacian with Robin boundary conditions. We summarise our results in the following
\begin{theorem}
	If $\f{R_\eps^N}{\eps^{N-1}}\to 1$ as $\eps\to 0$, then the solutions $u_\eps$ of \eqref{eq:graph_problem} satisfy $\big\|u_\eps-\V^\Gamma u\big\|_{L^2(\Om_\eps)}\to 0$, where $u\in H^1(\Gamma)$ solves
	\begin{align}\label{eq:borderline_graph_limit}
		\begin{cases}
			(-\Delta+z+\mu) u = f &\text{ on }\Gamma\\
			\quad\;\,\, \sum_{e\ni v} u'_e(v) = (z+\mu) \f{|V|}{|\Om_0|} u(v), &\text{ at each vertex }v
		\end{cases}
	\end{align}
	In particular, the strange term $\mu$ enters the vertex condition of the limit problem.
\end{theorem}
\section{Conclusion}
We have shown that the classical result by \cite{CM} also holds in a thin domain shrinking towards an interval or a graph. Furthermore, norm-resolvent convergence {holds} in the sense of Theorem \ref{normconv} and convergence of eigenvalues. Several generalisations {naturally arise}. First, the author believes that the norm convergence result generalises to unbounded domains (that is, when the limit domain is an unbounded interval). A suitable modification of the argument in \cite{CDR} or \cite{KP} seems like a reasonable approach. 

Second, the curious effect of the ``strange term'' $\mu$ appearing in the vertex condition observed in section \ref{sec:borderline} requires further study. Spectral convergence and abstract operator estimates will be the subject of future work.


\begin{thebibliography}{20}
	\bibitem[AP10]{AP10}
	    {\sc J. M. Arrieta, and M. C. Pereira},
	    {\em Elliptic problems in thin domains with highly oscillating boundaries}, SeMA J., 51(1) (2010), pp.~17--24.
	\bibitem[AV14]{AV14}
		{\sc J. M. Arrieta, and M. Villanueva-Pesqueira},
		{\em Locally periodic thin domains with varying period},
		C. R. Math., 352(5) (2014), pp.~397--403.
	\bibitem[AV16]{AV16}
		{\sc J. M. Arrieta, and M. Villanueva-Pesqueira},
		{\em Thin domains with non-smooth periodic oscillatory boundaries},
		J. Math. Anal. Appl., 446(1)  (2017), pp.~130--164.
	\bibitem[Boe17]{B17}
		{\sc S.~B\"{o}gli},
		{\em Convergence of sequences of linear operators and their spectra},
		Integral Equations Operator Theory, 88(4) (2017), pp.~559--599.
	\bibitem[Boe18]{B18}
		{\sc S.~B\"{o}gli},
		{\em Local convergence of spectra and pseudospectra},
		J. Spectr. Theory, 8(3) (2018), pp.~1051--1098.
	\bibitem[BCD16]{BCD}
		{\sc D.~Borisov, G.~Cardone and T.~Durante},
		{\em Homogenization and norm-resolvent convergence for elliptic operators
	 	 in a strip perforated along a curve},
		Proc. Roy. Soc. Edinburgh Sect. A, 146(6) (2016), pp.~1115--1158.
	\bibitem[CM97]{CM}
		{\sc D. Cioranescu, and F. Murat},
		{\em A Strange Term Coming From Nowhere},
		Progr. Nonlinear Differential Equations Appl., 31 (1997), pp.~45--93.
	\bibitem[CDR17]{CDR}
		{\sc K.~{Cherednichenko}, P.~{Dondl} and F.~{R{\"o}sler}},
		{\em Norm-Resolvent Convergence in Perforated Domains},
		Asymptot. Anal., 110(3--4) (2018), pp.~163--184.
	\bibitem[EP05]{EP}
		{\sc P.~Exner and O.~Post},
		{\em Convergence of spectra of graph-like thin manifolds},
		 J. Geom. Phys., 54(1) (2005), pp.~77--115.
	\bibitem[IOS89]{IOS}
		{\sc G. A. Iosif'yan, O. A. Oleinik, and A. S. Shamaev},
		{\em  Mathematical Problems in Elasticity and Homogenization},
		Elsevier Science, Netherlands, 1992.
	\bibitem[KP17]{KP}
  		{\sc A. Khrabustovskyi, and O. Post},
  		{\em Operator estimates for the crushed ice problem},
  		Asymptot. Anal., 110(3-4) (2018), pp.~137--161.
	\bibitem[KZ03]{KZ}
		{\sc P.~Kuchment, H.~Zeng},
		{\em Asymptotics of spectra of {N}eumann {L}aplacians in thin domains},
		in Advances in Differential Equations and Mathematical Physics: UAB International Conference, Differential Equations and Mathematical Physics, {B}irmingham, {AL}, 2002, pp.~199--213.
	\bibitem[MK64]{MK}
		{\sc V. A.  Marchenko and E. Ya. Khruslov},
  		{\em Boundary-value problems with fine-grained boundary [in Russian]},
  		Mat. Sb. (N.S.), 65(107):3 (1964), pp.~458--472.
  	\bibitem[MS10]{MS10}
		{\sc J. S. Mart\'\i n and L. Smaranda},
		{\em Asymptotics for eigenvalues of the {L}aplacian in higher
              dimensional periodically perforated domains},
        Z. Angew. Math. Phys., 61(3) (2010), pp.~401--424.
    \bibitem[MP10]{MP10}
		{\sc T. A. Mel'nyk and A. V. Popov},
		{\em Asymptotic analysis of boundary-value problems in thin perforated domains with rapidly varying thickness},
		Nonlinear Oscill., 13(1) (2010), pp.~57--84. 
	\bibitem[MP12]{MP12}
		{\sc T. A. Mel'nyk and A. V. Popov},
		{\em Asymptotic analysis of boundary value and spectral problems in thin perforated regions with rapidly changing thickness and different limiting dimensions},
		Sb. Math 203(8) (2012), pp.~1169--1195. 
	\bibitem[MNP13]{MNP}
		{\sc D.~Mugnolo, R.~Nittka and O.~Post},
		{\em  Norm convergence of sectorial operators on varying {H}ilbert spaces},
		Oper. Matrices, 7(4)  (2013), pp.~955--995.
	\bibitem[Naz10]{N10}
		{\sc S. A. ~Nazarov},
		{\em Opening of a Gap in the Continuous Spectrum of a Periodically Perturbed Waveguide},
		Math. Notes 87(5) (2010) pp.~738--756.
	\bibitem[Pas06]{Pas06}
		{\sc S.~E. Pastukhova},
		{\em Some estimates from homogenized elasticity problems},
		Dokl. Math., 73(1) (2006), pp.~102--106.
	\bibitem[Pos06]{P06}
		{\sc O. Post},
		{\em Spectral Convergence of Quasi-One-Dimensional Spaces}
		Ann. Henri Poincar\'{e} 7(5)  (2006), pp.~933--973.
	\bibitem[Pos12]{P12}
		{\sc O. Post},
		{\em Spectral analysis on graph-like spaces},
		Springer, Heidelberg (2012).
  	\bibitem[RT75]{RT}
  		{\sc J. Rauch, M. Taylor},
  		{\em Potential and scattering theory on wildly perturbed domains},
  		J. Funct. Anal., 18 (1975), pp.~27--59.
	\bibitem[Stu70]{Stummel1}
		{\sc F. Stummel},
		{\em Diskrete Konvergenz linearer Operatoren I},
		Math. Ann. 190 (1970), pp.~45--92.
	\bibitem[Stu72]{Stummel2}
		{\sc F. Stummel},
		{\em Diskrete Konvergenz linearer Operatoren II},
		Math. Z. 120 (1971), pp.~231--264.
	\bibitem[Vai81]{V81}
		{\sc G. M. Vainikko},
		{\em Regular convergence of operators and approximate solution of equations},
		J. Sov. Math. 15(6) (1981), pp.~675--705.
	\bibitem[Zhi00]{Jikov2000}
		{\sc V.~V. Zhikov},
		{\em On an extension and an application of the two-scale convergence method},
		Mat. Sb., 191(7) (2000), pp.~31--72.
\end{thebibliography}
\end{document}